%% file: preprint_elamg.tex
\pdfoutput=1
\documentclass{article}

\usepackage[numbers,sort&compress]{natbib}
\input{header.tex}
\input{header2.tex}
\usepackage[section]{placeins}

\newcommand{\kdim}{k}

\begin{document}

\title{An algebraic multigrid method based on an auxiliary topology with edge matrices}

\author{Lukas Kogler, Joachim Sch\"oberl \\ Institute of Analysis and Scientific Computing, Vienna Technical University, Austria}

\maketitle

\begin{abstract}
  This paper introduces a novel approach to algebraic multigrid methods for large systems of linear equations coming from finite element discretizations of certain elliptic second order partial differential equations.
Based on a discrete energy made up of edge and vertex contributions, we are able to develop coarsening criteria that guarantee two-level convergence even for systems of equations.
This energy also allows us to construct prolongations with prescribed sparsity pattern that still preserve kernel vectors exactly.
These allow for a straightforward optimization that simplifies parallelization and reduces communication on coarse levels.
Numerical experiments demonstrate efficiency and robustness of the method and scalability of the implementation.
\end{abstract}
\date{}

\section{Introduction}\label{sec:intro}
\input{intro.tex}

\section{Preliminaries}\label{sec:fem}
\input{fem.tex}

\section{Algebraic Multigrid Methods}\label{sec:amg}
\input{amg.tex}

\section{Discrete Energy}\label{sec:deng}
\input{deng.tex}

\section{Coarsening}\label{sec:crs}
\input{crs2.tex}

\section{Auxiliary Smoothed Prolongation}\label{sec:prol}
\input{prol.tex}

\section{Numerical Results}\label{sec:res}
\input{res.tex}

\section{Conclusions}\label{sec:conc}
We have introduced a new notion of discrete energy made up of vertex and edge contributions based on local penalization of deviation from
zero energy functions.
The existence of suitable edge matrices on the finest level is a result from \cite{CMIII}, and coarse grid edge- and vertex matrices have
been identified by us.
From this energy, we derived criteria for coarsening algorithms that are proven to be stable even for difficult elasticity problems.
A concrete algorithm that incorporates these criteria has been developed based on one from literature (\cite{Napov2012}).
Furthermore, the edge matrices making up the discrete energy can be harnessed to construct a form of smoothed prolongation the sparsity pattern of which
can be controlled in a natural way.
The scalability of the method and its implementation have been demonstrated by a series of numerical experiments.

\section{Acknowledgments}\label{sec:ack}
\noindent
The authors acknowledge support by the Austrian Science Fund (FWF) through the projects F65 and W1245.

\noindent
The computational results have been achieved using the Vienna Scientific Cluster (VSC).

\bibliography{preprint_elamg}%

\end{document}

%% file: header.tex
\usepackage{amsmath}
\usepackage{amssymb}
\usepackage{amsthm}
\usepackage{enumerate} 
\usepackage[noend]{algpseudocode}
\usepackage{algorithm,algorithmicx}
\usepackage{bm} 

\usepackage{tikz}

\newtheorem{theorem}{Theorem}[section]

\newtheorem{lemma}{Lemma}[section]

\theoremstyle{definition}
\newtheorem{definition}{Definition}[section]

\theoremstyle{remark}
\newtheorem{remark}{Remark}[section]

\theoremstyle{definition}
\newtheorem{problem}{Problem}[section]

\theoremstyle{remark}
\newtheorem{notation}{Notation}[section]

\usepackage{adjustbox}
\usepackage{calc}  

\usepackage[margin=2.5cm,footskip=0.25in]{geometry}

\providecommand{\keywords}[1]
{
  \small	
  \textbf{\textit{Keywords---}} #1
}

%% file: header2.tex
\DeclareMathOperator{\Div}{div}
\DeclareMathOperator{\Skew}{skew}
\DeclareMathOperator{\trace}{tr}

\usepackage{bm} 
\newcommand{\vof}[1]{\bm{#1}}
\newcommand{\mat}[1]{\bm{#1}}
\newcommand{\vu}{\vof{u}}
\newcommand{\vv}{\vof{v}}
\newcommand{\vw}{\vof{w}}
\newcommand{\vx}{\vof{x}}
\newcommand{\vb}{\vof{b}}

\newcommand{\vr}{\vof{r}}

\newcommand{\mA}{\mat{A}}
\newcommand{\mB}{\mat{B}}
\newcommand{\mI}{\mat{I}}

\newcommand{\mM}{\mat{M}}
\newcommand{\mD}{\mat{D}}
\newcommand{\mP}{\mat{P}}
\newcommand{\mE}{\mat{E}}
\newcommand{\mR}{\mat{R}}

\newcommand{\mhA}{\hat{\mat{A}}}
\newcommand{\mhM}{\hat{\mat{M}}}
\newcommand{\mhD}{\hat{\mat{D}}}

\newcommand{\mtA}{\tilde{\mat{A}}}
\newcommand{\mtD}{\tilde{\mat{D}}}
\newcommand{\mtM}{\tilde{\mat{M}}}

\newcommand{\mtP}{\tilde{\mat{P}}}

\newcommand{\mcV}{\mathcal{V}}
\newcommand{\mcC}{\mathcal{C}}
\newcommand{\mcD}{\mathcal{D}}
\newcommand{\mcN}{\mathcal{N}}
\newcommand{\mcS}{\mathcal{S}}
\newcommand{\mcT}{\mathcal{T}}

\newcommand*{\Let}[2]{\State #1 $\gets$ #2}
\newcommand{\Matrix}[1]{\left(\begin{matrix}#1\end{matrix}\right)}
\newcommand{\norm}[1]{\left\lVert#1\right\rVert}

\newcommand{\KTG}{K_{TG}}

%% file: intro.tex
Algebraic multigrid (AMG) methods are one of the most popular classes of methods for the iterative solution of large sparse systems of linear equations.
Besides the classical approach \cite{Brandt1982AlgebraicM, RS_1984, Ruge1987}, various other AMG methods have been developed.
For an overview of the mainstream methods, see e.g \cite{amg_ov_Xu2017}.
The ones that are most relevant to the discussion here are smoothed aggregation algebraic multigrid (SAAMG),
\cite{Vanek_SA_1992, Vanek_SA_1996, Vanek_SA_1999, Vanek_SA_2001},
and also algebraic multigrid based on computational molecules (AMGm), \cite{CMI, CMII, CMIII}.
The latter is an attempt to combine the classical approach with algebraic multigrid based on element interpolation (AMGe), \cite{Brezina_eAMG_2001}.
It expands on the idea of element preconditioning introduced in \cite{Sch_elpre_2001}, where an AMG method was applied not to the system matrix itself, but to an auxiliary matrix obtained from preconditioning of the element matrices.
AMGm is based on the concept of small edge matrices that can be assembled into an auxiliary matrix.
In contrast to AMGe, which is based on element matrices that need additional topological information on coarse levels, AMGm only needs topological information stored in the matrix itself.
In \cite{CMIII}, a coarsening strategy for linearized elasticity problems based on the solution of small eigenvalue problems obtained from these edge matrices was introduced.

In this paper we propose a method of smoothed aggregation type that is, just as AMGm, based on the concept of edge matrices.
The method is suited for both scalar $H^1$ and linearized elasticity problems, and, in fact, on the theoretical level, hardly makes a distinction between them at all.
In the smoothed aggregation method a tentative coarse space is defined by restrictions of kernel vectors to agglomerates.
For linear elasticity, this means that coarse spaces contain both displacement and rotational degrees of freedom.
The introduction of edge matrices that explicitly include these rotational degrees of freedom is the crucial ingredient for the method proposed here.
We show that these edge matrices can be used to define both a robust coarsening algorithm based on edge-wise eigenvalue problems and a new variation of the smoothed prolongation used in the smoothed aggregation method.
This auxiliary smoothed prolongation can be used to control the sparsity of the prolongation operator in a natural fashion, which allows for a series of optimizations.
The validity of the approach and efficiency of the implementation are shown through computations.

Besides the works on AMGm, there is also a strong relation to \cite{Napov2010} and \cite{Napov2012}, where, based on a theoretical result from \cite{Falgout2004}, an aggregation algorithm for scalar problems that is based on similar eigenvalue problems was introduced.
In order to formulate the eigenvalue problems in that method, a particular nonnegative splitting of the global matrix is required.
While finding such a splitting of the actual system matrix is relatively straightforward only for scalar problems, there is a canonical splitting of the auxiliary matrix introduced here.
In fact, the eigenvalue problems in the coarsening algorithm proposed here are similar to ones that appear in the algorithm from \cite{Napov2012} when applied to the auxiliary matrix with this splitting, however they feature an additional stabilization that turns out to be essential for elasticity problems.
In that respect, they are more reminiscent of a coarsening criterion used in AMGm that is based on computational molecules, introduced in \cite{CMIII}.

\paragraph{Outline}
After stating the equations and discretizations we consider in section \ref{sec:fem}, section \ref{sec:amg} will discuss aggregation based AMG methods as well as a formulation of the tentative prolongation used in these methods that is the basis for the definition of edge matrices in section \ref{sec:deng}.
Based on these, a coarsening algorithm will be developed in section \ref{sec:crs} and the aforementioned auxiliary smoothed prolongation in section \ref{sec:prol}.
Finally, numerical results are presented in section \ref{sec:res}.

%% file: fem.tex
We briefly summarize the considered equations and the chosen discretizations, see \cite{Braess2007} for a more detailed description of finite element methods in general.
From now on, the dimension $d$ will be $2$ or $3$ and the domain $\Omega\subset\mathbb{R}^{d}$ will have the Dirichlet boundary $\Gamma_D\subseteq\partial\Omega$.
On the continuous level, we consider linear equations of the form:
\begin{problem}\label{prob:cont}
  With $V = H^1(\Omega) \text{ or }V = [H^1(\Omega)]^d$, an elliptic and continuous bilinear form $a(\cdot, \cdot)$, and a continuous linear form $g(\cdot)$ on $V$,
  find $u\in V_0 = \{u \in V : u_{|\Gamma_D}=0\}$, such that
  \begin{alignat}{1}\label{eq:lineq}
    a(u, v) = g(v) \quad\quad\forall v\in V_0.
  \end{alignat}
\end{problem}

On the discrete level we only consider discretization with the standard, lowest order, conforming finite element method.
For $\mathcal{T}$, a shape-regular triangulation of $\Omega$ that resolves $\Gamma_D$, and its set of vertices, $\mathcal{V} = \{v_1,\ldots, v_{n}\}$, in the scalar case that is
\begin{alignat}{2}\label{def:vhscal}
  V_h &:= \left\{u \in H^1(\Omega) : u_{|T} \in P_1(T),~\forall T\in\mathcal{T}\right\},
\end{alignat}
where  $P_1(T)$ is the set of polynomials of order 1 on $T$.
This space is spanned by the standard hat functions associated to vertices
\begin{alignat}{2}\label{def:phiscal}
  \{ \varphi^i : i\in\{1,\ldots, n\}\} \quad\text{where}\quad \varphi^i_{|T}\in P_1(T)~\forall T\in\mathcal{T},\varphi^i(v_j) = \delta_{ij}\quad\forall v_j\in\mathcal{V}.
\end{alignat}
For the vector case,
\begin{alignat}{2}\label{def:vhvec}
  V_h &:= \left\{u \in [H^1(\Omega)]^d : u_{|T} \in [P_1(T)]^d,~\forall T\in\mathcal{T}\right\}
\end{alignat}
is spanned by basis functions that are just the scalar ones times unit vectors $\vec{e}_l$,
\begin{alignat}{2}\label{def:phivec}
  \{ \varphi^{i,j} : i\in\{1\ldots n\}, j\in\{1\ldots d\}\}\quad\text{where}\quad \varphi^{i,l}_{|T}\in [P_1(T)]^d~\forall T\in\mathcal{T},\varphi^{i,l}(v_j) = \delta_{ij}\vec{e}_l\quad\forall v_j\in\mathcal{V} .
\end{alignat}
By Galerkin-projection of (\ref{eq:lineq}) onto $V_{h, 0} := V_h\cap V_0$, we obtain the system of linear equations that is to be solved,
\begin{align}\label{eq:auf}
  \mA\vu = \vb
\end{align}
where $\mA_{ij} = a(\varphi^i, \varphi^j)$ and $\vb_i = g(\varphi^i)$ for the basis functions $\varphi^i$ of $V_{h, 0}$.
As $a(\cdot, \cdot)$ is symmetric and elliptic, $\mA$ is a symmetric positive definite matrix.
The problem we consider for the scalar case is:
\begin{problem}\label{prob:gg}
  Solve the system of linear equations (\ref{eq:auf}) with $V_h$ as in (\ref{def:vhscal}) and
  \begin{alignat}{1}\label{eq:gg}
    a(u, v) &= \int_{\Omega}\alpha~\nabla u\cdot\nabla v + \beta~uv~dx
  \end{alignat}
  Here, $\alpha>0$ is bounded away from zero and either $|\Gamma_D|>0$ or $\beta>0$ (in some part of $\Omega$).
\end{problem}
In the vector valued case we are interested in linearized elasticity:
\begin{problem}\label{prob:epseps}
  Solve the system of linear equations (\ref{eq:auf}) with $V_h$ as in (\ref{def:vhvec}) and
  \begin{alignat}{1}\label{eq:epseps}
    a(u, v) = \int_{\Omega}\mu~\varepsilon(u):\varepsilon(v) + \lambda~\Div(u)\Div(v) + \beta~u\cdot v ~dx.
  \end{alignat}
  with $\varepsilon(u)=\frac{1}{2}(\nabla u + \nabla u^T)$.
  Now $\mu>0$ is bounded away from zero and $\lambda\geq 0$ and $\beta\geq 0$.
  Additionally, we exclude the almost incompressible case, where $\mu\ll\lambda$, where, as is well known, the problem becomes extremely ill-conditioned (see \cite{Braess2007}).
\end{problem}

For the coordinate vector of a function $u$ or $u_h$ in $V_h$ with respect to the basis functions given in (\ref{def:phiscal}) or (\ref{def:phivec}) above, we write $\vof{u}$.
Owing to the clear association of basis functions and vertices, we will, for the elasticity problem, consider $\mA$ not as an element of $\mathbb{R}^{nd\times nd}$ but rather of $(\mathbb{R}^{d\times d})^{n\times n}$, so an ``entry'' $\mA_{ij}\in\mathbb{R}^{d\times d}$ stands for the $d\times d$ sub-matrix of $\mA$ corresponding to basis functions associated to vertex $v_i$ and $v_j$.
Therefore, by saying that a matrix is diagonal, we mean it to be block-diagonal with respect to vertices.
Similarly, vectors $\vof{u}\in\mathbb{R}^{nd}$ are interpreted as elements of $(\mathbb{R}^{d})^{n}$, so an ``entry'' $\vof{u}_i\in\mathbb{R}^d$ is the sub-vector of $\vof{u}$ consisting of entries associated to vertex $v_i$.
There will also be matrices and vectors (on coarse levels) that are to be understood as having a bigger "block size" than $d$, which will be clear from context. 
For the most part, vectors and matrices that are to be understood this way are written bold, $\mA, \vx$ to differentiate them from generic ones, $A, x$.
For symmetric matrices $A$ and $B$, we write $A\leq B$ if $\forall x: x^TAx \leq x^TBx$ and $A<B$ defined in the same manner.

%% file: amg.tex
Multigrid methods for the solution of (\ref{eq:auf}) have two components: smoothers, some sort of computationally cheap iterative methods, and prolongations.
A prolongation matrix is the representation of an embedding operator of some coarse space of smaller dimension into the original, fine space.
A simple two-grid method (Algorithm \ref{alg:tg}) consists of one application of a smoother, restriction of the residual to the coarse space, the exact solution of the projected sytem on the coarse space, embedding of the coarse solution into the fine space and finally application of the transposed smoother which makes the whole iteration symmetric.
For a thorough discussion of multilevel methods in general, see \cite{mlfvass}.

\begin{algorithm}
  \begin{algorithmic}
    \Let{$\vx$}{$\vx + \mM^{-1}(\vb - \mA\vx)$} \Comment{pre-smoothing}
    \Let{$\vx$}{$\vx + \mP(\mP^T\mA\mP)^{-1}\mP^T(\vb-\mA\vx)$} \Comment{coarse grid correction}
    \Let{$\vx$}{$\vx + \mM^{-T}(\vb - \mA\vx)$} \Comment{post-smoothing}
  \end{algorithmic}
  \caption{Two-grid method}
  \label{alg:tg}
\end{algorithm}

The iteration matrix of the two-grid method is
$$
\mE_{TG} = \mI - \mB_{TG}^{-1}\mA = (\mI-\mM^{-T}\mA)(\mI-\mP(\mP^T\mA\mP)^{-1}\mP^T\mA)(\mI-\mM^{-1}\mA)
$$
where
$$
\mB_{TG}^{-1} = \widetilde{\mM}^{-1} + (\mI-\mA\mM^{-T})\mP(\mP^T\mA\mP)^{-1}\mP^{T}(\mI-\mM^{-1}\mA)
$$
with the symmetrized smoother
$$
\widetilde{\mM} := \mM(\mM+\mM^T-\mA)^{-1}\mM^T.
$$
For algebraic multigrid methods, a simple smoother $\mM$ is used - (Block-)Gauss-Seidel type and polynomial smoothers are popular.
Then, clearly, the key ingredient is the prolongation matrix $\mP$, for which it is of great importance that its range spans the kernel vectors of $\mA$ (without boundary conditions, and disregarding $L^2$ terms).
In aggregation-based methods, a first, tentative prolongation is constructed by finding some partition of the fine degrees of freedoms into so-called agglomerates, and then defining coarse degrees of freedom by the restrictions of kernel vectors to these agglomerates.
In order to motivate the definition of edge matrices in the next section, we give a somewhat roundabout formulation of this tentative, or unsmoothed, prolongation.
For that, suppose that  for any two vertices $v_i$ and $v_j$ there exists a matrix $Q^{v_i\rightarrow v_j}\in\mathbb{R}^{k\times k}$, where $k$ is the number of degrees of freedom per vertex, such that
$$
\mA\vu = 0 \quad\Leftrightarrow\quad Q^{v_i\rightarrow v_j}\vu_i = \vu_j \quad\forall v_i, v_j\in\mathcal{V}.
$$
When a global kernel vector can be identified by its entries at a single vertex, $Q^{v_i\rightarrow v_j}$ represents the transform of coordinates associated to one vertex into those an energy free extension would have at another vertex.
For example, for the scalar case (Problem \ref{prob:gg}), take $Q^{v_i\rightarrow v_j}\in\mathbb{R}^{1\times 1} = 1$, and the above expression is just
$$
\mA\vu = 0 \quad\Leftrightarrow\quad \vu_i = \vu_j \quad\forall v_i, v_j\in\mathcal{V}.
$$
This obviously captures the constant vector, which is the single kernel vector of $\mA$.
The situation is less clear for the equations of linearized elasticity, Problem \ref{prob:epseps}, where, in dimension $d$, we have $d$ degrees of freedom per vertex, but $k = d+\frac{d(d-1)}{2}$ rigid body modes that form the kernel.
With
$$
\Skew(\vec{t}) := \Matrix{-t_1\\t_0}
$$
in two dimensions, and
$$
\Skew(\vec{t}) := 
\Matrix{
  0    & -t_3 & t_2 \\
  t_3  & 0    & -t_1 \\
  -t_2 & t_1  & 0 \\
}
$$
in three dimensions, these are the functions
$$
u(\vec{x}) = \vec{d} + \Skew(\vec{x})~\vec{r}, \quad \vec{d}\in\mathbb{R}^d, \vec{r}\in\mathbb{R}^{k-d}.
$$
For this reason, let us take a detour and suppose we were dealing with some discretization that has exactly $k$ degrees of freedom per vertex, which are nodal in the sense that
$$
u(\vec{x}) = \vec{d} + \Skew(\vec{x})~\vec{r} \quad\Leftrightarrow\quad \vu_i = \Matrix{ u(\vec{x}_i) \\ \vec{r}  }\quad\forall v_i\in\mcV.
$$
where $x_i$ is the position of vertex $v_i$.
In this scenario, every vertex would have associated rotational degrees of freedom in addition to displacement degrees of freedom. Then, from
$$
u(\vec{y}) = \vec{d} + \Skew(\vec{y})~\vec{r} = u(\vec{x}) + \Skew(\vec{y}-\vec{x})~\vec{r}
$$
it follows that
\begin{align}
\mA\vu = 0 \quad\Leftrightarrow\quad
Q^{v_i\rightarrow v_j}\vu_i = \Matrix{
  I & \Skew{\vec{t}^{ij}} \\
  0 & I
}
\vu_i = \vu_j \quad\forall v_i, v_j\in\mathcal{V}\label{eq:qij_epseps}
\end{align}
where $\vec{t}^{ij} = \vec{x}_j - \vec{x}_i$ is the vector connecting vertices $v_i$ and $v_j$.
\begin{notation}
  The number of degrees of freedom per vertex will be referred to as $k$. It should be clear from context whether $k=d$ or $k=d+\frac{d(d-1)}{2}$.
\end{notation}
Expressing the change in coordinates from one vertex to another by such matrices, we can formulate a tentative prolongation, as used in aggregation based AMG methods.
\begin{definition}[tentative prolongation]
  \label{def:pwp}
  For the set of fine vertices $\mathcal{V} = \{v_i, i\in\mathcal{I}_f \}$, a set $\mathcal{D}\subseteq\mathcal{V}$
  and a partition $\mathcal{C} = \{C_I, I\in \mathcal{I}_c\}$ of $\mathcal{V}\setminus\mathcal{D}$ such that $C_I\cap C_J = \emptyset, I\neq J$,
  define a coarse vertex $v_I$ at some position $x_I$ for every $C_I\in\mathcal{C}$.
  Then the tentative prolongation $\mP\in(\mathbb{R}^{k\times k})^{|\mathcal{I}_f|\times|\mathcal{I}_c|}$ is defined by
  \begin{align}
    \mP_{iJ} :=
    \begin{cases}
      Q^{v_J\rightarrow v_i} & \quad \text{if } i\in C_J, \\
      0 & \quad \text{else}
    \end{cases}. \label{eq:pwp}
  \end{align}
\end{definition}

\begin{notation}
In order to more easily differentiate between fine and coarse level, we use lower case letters for indices associated with the fine level and upper case letters for those associated with the coarse level.
\end{notation}

\begin{remark}
  Usually, aggregation-based algebraic multigrid methods orthogonalize the restrictions of kernel vectors to agglomerates.
  This has a similar effect to moving $x_I$, which amounts to choosing a central point with respect to which coarse grid rotations are defined.
  Therefore, it is reasonable to choose coarse vertices somehow towards the "center" of agglomerates, e.g. $x_I = \frac{1}{|C_I|}\sum_{i\in C_I}x_i$.
\end{remark}
When using discretization (\ref{def:vhvec}) for linearized elasticity, there are no rotational degrees of freedom on the finest level, so Definition \ref{def:pwp}
cannot be used immediately.
However, we can express a prolongation from a coarse level that does include rotational degrees of freedom by simply restricting the above one to rows
corresponding to displacement degrees of freedom:
\begin{definition}[finest level elasticity prolongation]
  \label{def:pwp_el}
  Let $\mP\in\left(\mathbb{R}^{k\times k}\right)^{|\mathcal{V}|\times |\mathcal{C}|}$ be a tentative prolongation for Problem \ref{prob:epseps} as in Definition \ref{def:pwp} with $Q^{v_I\rightarrow v_j}$ defined as in (\ref{eq:qij_epseps}),
  $k=d+\frac{d(d-1)}{2}$ and $\mat{R}\in\left(\mathbb{R}^{d\times k}\right)^{|\mathcal{V}|\times |\mathcal{V}|}$ with $\mat{R}_{ij}=\delta_{ij}\Matrix{I_d & \mat{0}}$ with the $d\times d$ identity matrix $I_d$.
  Then the finest level tentative prolongation $\mP^{f}\in(\mathbb{R}^{d\times k})^{|\mathcal{V}|\times|\mathcal{C}|}$ for Problem \ref{prob:epseps} is given by
  \begin{align}
    \mP^f := \mat{R}\mP. \label{eq:pwp_fl}
  \end{align}
\end{definition}

\input{sing2.tex}

%% file: sing2.tex
A potential issue with (\ref{eq:pwp_fl}) is that without further restrictions to $\mcC$, $\mat{R}\mP$ does not necessarily have full column rank.
This happens in three dimensions when all vertices in an agglomerate lie on a straight line, and in both two and three dimensions when an agglomerate
consists of only one vertex.
This means that coarse level matrices can be singular if no extra care is taken during the coarsening process to make sure that none of the cases above occur.
As we aim to provide a general purpose coarsening algorithm in this paper, we consider an alternative solution to this problem.
In principle, one could modify the prolongations themselfs - for a series of prolongations $\mP_0,\ldots, \mP_{L-1}$, let, for $l>0$, the columns of the matrix $\mE_l$ consist of an orthonormal base of $\ker{\left(\mR\mP_o\cdot\mP_1\cdot\cdot\cdot\mP_{l-1}\right)}^{\perp}$ and set $\mE_0:=\mR^T$.
Then the "restricted" prolongations $\mtP_l := \mE_{l-1}^T\mP_{l-1}\mE_l$ induce regular coarse grid matrices $\mtA_0:=\mA,~\mtA_l:=\mtP_{l-1}^T\mtA_{l-1}\mtP_{l-1}$.
However, the number of degrees of freedom per vertex is not necessarily constant anymore.
Lemma \ref{thm:eet_is_diag} characterizes the kernel of the coarse grid matrices and shows that, unsurprisingly, the matrices $\mE_l$ are feasible to compute.
\begin{lemma}\label{thm:eet_is_diag}
  Let $\mat{A}_0$ be the finite element matrix for Problem \ref{prob:epseps}, and $\mA_1,\ldots, \mA_L$ be a series of coarse grid matrices,
  obtained by Galerkin projection, $\mA_1$ by $\mR\mP_0$ as in \ref{def:pwp_el}, and $\mA_2,\ldots, \mA_L$ by $\mP_1,\ldots, \mP_{L-1}$ from
  Definition \ref{def:pwp}. 
  Then there exists an orthonormal base of $\ker{\mat{A}}_l^{\perp}$ that consists of vectors that have nonzero entries only associated to a single vertex each.
  The matrix $\mat{E}_l$, the columns of which consist of that basis, where all vectors belonging to a vertex are numbered consecutively, has the block structure
  \[
    \mE_l =
    \Matrix{ \mE_{l,11} & & & \\
      & \mE_{l,22} & & \\
      & & \ddots & \\
      & & & \mE_{l,mm}\\
    },
  \]
  where $m=|\mcV_l|$ is the number of vertices on level $l$ and $\mE_{l, ii}\in\mathbb{R}^{\kdim\times n_i}$ with $n_i\leq\kdim$.
\end{lemma}
\begin{proof}
  Firstly, $\mA_l = \mP_{l-1}^T\mA_{l-1}\mP_{l-1} = \ldots = \mP_{l-1}^T\cdots\mP_{0}^T\mR^T\mA\mR\mP_0\cdots\mP_{l-1} = \mtP^T\mR^T\mA\mR\mtP$,
  where $\mtP:=\mP_{0}\cdot\mP_1\cdots\mP_{l-1}$ is itself a prolongation of the type of Definition \ref{def:pwp}.
  It is therefore enough to consider only the case where there is only one coarse grid matrix $\mat{A}_1 = (\mat{R}\mat{P})^T\mat{A}\mat{R}\mat{P}$.
  To proof the statement, we need to find a basis of $\ker{\mat{A}_1}=\ker{\mR\mP}$ where every basis vector only has entries associated to a single vertex.
  As $\mat{R}\mat{P}$ only has (at most) one column entry per row, for any vector $\vof{u}\in\ker{(\mat{R}\mat{P})}^\perp$ also its restriction to a single vertex is in the kernel.
  For any such base of $\ker{\mR\mP}$, therefore, also a maximal linearly independent subset of the restrictions of these vectors to single vertices form a base.
  Orthogonalization and normalization threreof gives a basis as is claimed to exist in the lemma. The block structure of $\mE$ follows immediately, where $n_i$ is the number of vectors in the basis associated to vertex $v_i$.
\end{proof}
Therefore, $\mE_l$ can be computed by finding the kernel of diagonal entries of $\mA_l$.
The disadvantage of this approach is that due to the varying number of degrees of freedom per vertex on coarse levels, we cannot make use of especially
efficient sparse matrix formats that store entries in blocks of fixed, static size.
We propose to use it for theory only, and to deal with singular coarse grid matrices in a different way.
Given some convergent smoother $\mtM_l$ for $\mtA_l$, we can define a "pseudo smoother" $\mM_l$ for $\mA_l$ as $\mM_l := \left(\mE_l\mtM_l\mE_l^T\right)$ where we make use of the pseudo inverse $(\cdot)^{\dagger}$ to define the iteration
$$
\vx\rightarrow \vx + \left(\mE\mtM_l\mE^T\right)^{\dagger}(\vb-\mA_l\vx) = \vx + \mE\mtM_l^{-1}\mE^T(\vb-\mA_l\vx).
$$
This can be seen as merely \textit{an implementation} of the smoothing iteration on the restricted matrix $\mtA_l$ using $\mtM_l$.
The advantage is that when $\mtM_l$ is the Jacobi or Gauss-Seidel smoother, $\mM_l$ is the Jacobi or Gauss-Seidel smoother for $\mA_l$, where
diagonals are pseudo inverted instead of inverted.
On the coarsest level, where we have to invert $\mtA_l$, observe that
$$
\left(\mE_l\mtA_l\mE_l + (\mI-\mE_l\mE_l^T) \right)^{-1} = \mE_l\mtA_L^{-1}\mE_l^T + (\mI-\mE_l\mE_l^T)
$$
and therefore
$$
\left(\mE_{l-1}\mP_{l-1}\right)\left(\mE_l\mtA_l\mE_l + (\mI-\mE_l\mE_l^T) \right)^{-1}\left(\mP_{l-1}^T\mE_{l-1}^T\right) = \mtP_{l-1}\mtA_L^{-1}\mtP_{l-1}.
$$
It is therefore justified to regularize the coarsest level matrix with the term $\mI-\mE_l\mE_l^T$, which, according to Lemma \ref{thm:eet_is_diag}, is feasible to compute
and, as seen above, does not disturb the coarse grid correction.
This allows us to work with the aforementioned sparse matrices with static block size, which should be seen as an implementation of the AMG method defined by the
restricted prolongations $\mtP_l$.

%% file: deng.tex
The idea of constructing an approximate splitting of the matrix $\mA$ into edge contributions has already been explored in
\cite{CMI, CMII, CMIII}. There, edge matrices of size $(2d)\times(2d)$ were used.
These edge matrices have to preserve the kernel of the global differential operator.
Using the explicit representation of these kernel vectors from the previous section, we can define edge matrices in an arguably simpler and more intuitive way:

\begin{definition}[auxiliary matrix and discrete energy]\label{def:deng}
  For $v_i, v_j\in\mathcal{V}$, let $Q^{v_i\rightarrow v_j}\in \mathbb{R}^{\kdim\times\kdim}$ be invertible such that
  $Q^{ij} := \sqrt{Q^{v_i\rightarrow v_j}}$ exists and $Q^{v_i\rightarrow v_j}Q^{v_j \rightarrow v_l} = Q^{v_i\rightarrow v_l}~\forall v_l\in\mathcal{V}$.
  For $v_i, v_j\in\mcV$, let $E^{ij}=E^{ji}\in \mathbb{R}^{\kdim\times\kdim}$ and $M^i\in \mathbb{R}^{\kdim\times\kdim}$ be symmetric and nonnegative.
  The auxiliary matrix $\mhA\in(\mathbb{R}^{\kdim\times\kdim})^{|\mathcal{V}|\times|\mathcal{V}|}$ defined by these edge and vertex matrices is:
  \begin{align}\label{eq:auxa}
    \mhA_{ij} = \begin{cases}
      M^{i} + \sum_{j \neq i} Q^{ij, T}E^{ij}Q^{ij}& i = j \\
      -Q^{ij, T}E^{ij}Q^{ji}& i \neq j 
      \end{cases}
  \end{align}
  $\mhA$ induces the discrete energy:
  \begin{align}\label{eq:deng}
    \norm{\vu}_{\mhA} = \sqrt{\sum_{v_i\in\mcV} \norm{M^i\vu_i}^2 + \sum_{v_i\neq v_j\in\mcV} \norm{Q^{ij}\vu_i - Q^{ji}\vu_j}_{E^{ij}}^2} .
  \end{align}
\end{definition}

We call $E^{ij}$ edge weight and $M^i$ vertex weight matrices. 
While the second term on the right in (\ref{eq:deng}) captures the energy induced by the differential operator ($e.g \norm{\varepsilon(u)}_{L_2}^2$), the first term is supposed to represent the $L_2$-term in (\ref{eq:gg}), (\ref{eq:epseps}) when $\beta>0$, or, as discussed shortly, naturally arises on coarse levels under the presence of Dirichlet boundary conditions.

\begin{remark}
  Although the use of $Q^{ij} = \sqrt{Q^{v_i\rightarrow v_j}}$ may seem unnecessarily complicated for the case of linearized elasticity,
  it is readily seen that
  \[
    Q^{ij} =
    \left(\begin{matrix}
        I & \frac{1}{2}\Skew(t^{ij})\\
        0 & I
      \end{matrix}\right).
  \]
  $Q^{ij}$ simply represents the transformation of coordinates given w.r.t. $v_i$ to ones given w.r.t. the midway point between $v_i$ and $v_j$.
  This definition imparts symmetry to expression (\ref{eq:deng}).
\end{remark}

For Problem \ref{prob:gg}, where $\kdim=1$ and $Q^{ij}=1$, under the condition that $\mA$ is an M-matrix, the choice of $E^{ij}=|\mA_{ij}|$, and $M^i=\mA_{ii} - \sum_{j\neq i} |\mA_{ij}| \geq 0$ gives an exact splitting of $\mA$ into edge- and vertex contributions.
Now, $\mA$ is not always an M-matrix, however, in practice, there seem to be no issues using this choice for $E^{ij}$ anyways and set $M^i=\max\{0, \mA_{ii} - \sum_j |\mA_{ij}|\}$.
Either way, the edge contributions simplify to
\begin{align}
  \norm{Q^{ij}\vu_i - Q^{ji}\vu_j}_{E^{ij}}^2 =
  E^{ij}
\left(\begin{matrix}
    \vu_i \\
    \vu_j
  \end{matrix}\right)^T
\left(\begin{matrix}
    1  & -1 \\
    -1 & 1
  \end{matrix}\right)
\left(\begin{matrix}
    \vu_i \\
    \vu_j
  \end{matrix}\right)
=
  E^{ij}(\vu_i-\vu_j)^2. \label{eq:gg_ectrb}
\end{align}
For Problem \ref{prob:epseps} in three dimensions, where $k=6$, when splitting displacement and rotational degrees of freedom, an edge contribution is
\begin{align}
  \norm{Q^{ij}\Matrix{\vu_i\\\vr_i} - Q^{ji}\Matrix{\vu_j\\\vr_j}}_{E^{ij}}^2
  =
  \norm{\begin{matrix}
        (\vu_i-\vu_j) - \frac{1}{2}\vec{t}^{ij}\times(\vr_i+\vr_j) \\
        \vr_i - \vr_j
      \end{matrix}}_{E^{ij}}^2
  =
  \norm{\begin{matrix}
      [\vu] - \vec{t}^{ij}\times\{\vr\} \\
      [\vr]
      \end{matrix}}_{E^{ij}}^2, \label{eq:epseps_ectrb}
\end{align}
where $[\cdot]$ denotes the jump and $\{\cdot\}$ the average along an edge.
\begin{figure}
  \begin{center}
    \begin{tabular}{cc}
      \includegraphics[scale=.8]{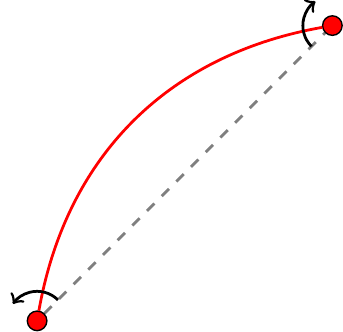} & \includegraphics[scale=.8]{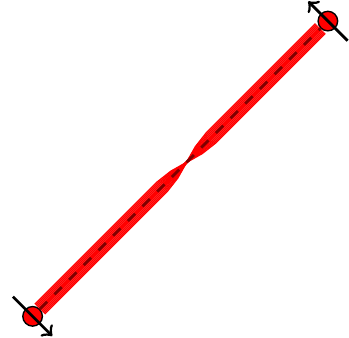} \\
      \text{"bending"} & \text{"twisting"} \\
      \includegraphics[scale=.8]{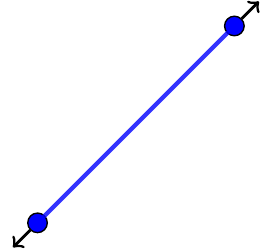} & \includegraphics[scale=.8]{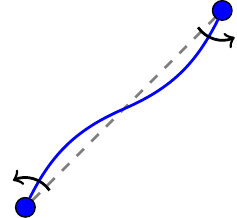} \\
      \text{"stretching"} & \text{"wiggling"}
    \end{tabular}
  \end{center}
  \caption{visualization of the components of the last term in (\ref{eq:epseps_ectrb})}
  \label{fig:enrg}
\end{figure}
The second entry of the above vector, $[\vr]$, stands for a ``bending'', or ``twisting'' along an edge, while the first, $[\vu] - \vec{t}^{ij}\times\{\vr\}$, describes compression/expansion (it's tangential part) and a kind of ``wiggling'' in it's normal part (see Figure \ref{fig:enrg})
\begin{remark}
  An energy functional made up of contributions similar to (\ref{eq:gg_ectrb}) is proposed in \cite{ef_eamg} as a component of the construction of the interpolation in the element-free AMGe method.
  An energy functional based on contributions of the form (\ref{eq:epseps_ectrb}) might be a good candidate with which to extend that approach also to elasticity problems.
\end{remark}
On the finest level, (where $k=d$), when $\beta=0$, Theorem 3.8 and the following Corollary 3.11 in \cite{CMII} show that an approximate splitting of $\mA$ into a sum of edge matrices of the form
\begin{align}
  \tilde{E}^{ij} = c^{ij} \Matrix{\vec{t}^{ij}\\-\vec{t}^{ij}}\Matrix{\vec{t}^{ij}\\-\vec{t}^{ij}}^T \label{eq:amge_ela_emats}
\end{align}
exists, with positive constants $c^{ij}\in\mathbb{R}$ and again the vector $\vec{t}^{ij} = \vec{x}_j - \vec{x}_i$ connecting vertices $v_i$ and $v_j$.
Theorem \ref{thm:deng_elast} shows that there is a choice of $E^{ij}$ in Definition \ref{def:deng} that reproduces these edge matrices, and thus, that in this case, the resulting auxiliary matrix $\mhA$ is equivalent to $\mA$.
\begin{theorem}[discrete energy for elasticity] \label{thm:deng_elast}
  Let $\mat{A}$ be the finite element matrix for Problem \ref{prob:epseps} with $\beta=0$, let $Q^{ij}$ be defined as in (\ref{eq:qij_epseps}) and
  let $\mat{E}\in(\mathbb{R}^{d\times k})^{|\mathcal{V}|\times |\mathcal{V}|}$ with
  \[
    \mat{E}_{ij} = \delta_{ij} \left(\begin{matrix} I_d \\ 0\end{matrix}\right).
  \]
  Then, there exist $c^{ij}\geq 0$ such that the choice of
  \[
    E^{ij} = c^{ij}
    \left(\begin{matrix}
        \vec{t}^{ij}\otimes \vec{t}^{ij}& 0\\
        0 & 0
      \end{matrix}\right)
  \quad\text{ and }\quad
    M^i = 0
  \]
  in Definition \ref{def:deng} defines an auxiliary matrix $\hat{\mat{A}}$ such that
  \begin{align}
    \norm{\vof{u}}_{\mat{A}} \sim \norm{\mat{E}\vof{u}}_{\mhA} = \sqrt{\sum_{v_i\neq v_j\in\mcV} \norm{Q^{ij}\Matrix{\vof{u}_i\\\vof{0}} - Q^{ji}\Matrix{\vof{u}_j\\\vof{0}}}_{E^{ij}}^2} \label{eq:fl_ela_equiv}
  \end{align}
  with constants depending only on the mesh regularity.
  Alternatively, setting $Q^{ij}=I_d\in\mathbb{R}^{d\times d}$ and $E^{ij} = c^{ij}\vec{t}^{ij}\otimes \vec{t}^{ij}$ in Definition \ref{def:deng} yields an auxiliary matrix $\tilde{\mat{A}}$ such that
  \begin{align}
    \norm{\vof{u}}_{\mat{A}} \sim \norm{\vof{u}}_{\tilde{\mat{A}}} = \sqrt{\sum_{v_i\neq v_j\in\mcV} \norm{\vof{u}_i- \vof{u}_j}_{E^{ij}}^2}. \label{eq:fl_ela_equiv_short}
  \end{align}
\end{theorem}

\begin{proof}
  As shown in Theorem 3.8 in \cite{CMII} for element matrices, and Corollary 3.11 for $\mA$, there exists an approximate splitting of $\mA$ into edge matrices $\tilde{\mat{E}}^{ij}$ of the form (\ref{eq:amge_ela_emats}) such that
  \begin{align}
    \norm{u}_{A} \sim \sqrt{\sum_{v_i\neq v_j\in\mcV} \norm{\Matrix{u_i\\u_j}}_{\tilde{E}^{ij}}^2}. \label{eq:amge_ela_eqiv}
  \end{align}
  It remains to show that any edge contribution in (\ref{eq:amge_ela_eqiv}) can be represented by one from (\ref{eq:fl_ela_equiv}) for some $E^{ij}$. Indeed, the choice of $E^{ij}$ given in the theorem, using the same $c^{ij}$, yields
  \[
    \norm{Q^{ij}\Matrix{u_i\\0} - Q^{ji}\Matrix{u_j\\0}}_{E^{ij}}^2
    =
    \Matrix{ u_i \\ u_j }^T
    c^{ij}\Matrix{
      \vec{t}^{ij}\otimes \vec{t}^{ij,T} & -\vec{t}^{ij}\otimes \vec{t}^{ij,T} \\
      -\vec{t}^{ij}\otimes \vec{t}^{ij,T} & \vec{t}^{ij}\otimes \vec{t}^{ij,T}
    }
    \Matrix{ u_i \\ u_j }
    =
    \norm{\Matrix{ u_i \\ u_j }}^2_{\tilde{E}^{ij}}.
  \]
  The second statement of the theorem follows from
  $
  \mE^T\mhA\mE = \mtA.
  $
\end{proof}


\begin{remark}
  For the problems considered in this paper, the rather crude choice of $c^{ij} = \frac{1}{d^2}\sum_{l,m\in\{0,\ldots, d\}}|(\mA_{ij})_{lm}|$ suffices.
  A more rigorous choice is used in \cite{CMIII}, where $c^{ij}$ are computed from off-diagonal entries of the finite element matrix $\mA$.
\end{remark}

So far, Theorem \ref{thm:deng_elast} is no real improvement over the results from \cite{CMII}.
The edge matrices $E^{ij}$ are rank 1 matrices, and only the upper left block of the transformation matrices $Q^{ij}$ is used.
However, combining this with Definition \ref{def:pwp} leads to a natural way of finding a splitting into edge matrices for the coarse grid matrix.

\begin{theorem}\label{thm:cemats}
  Let the fine level matrix $\mA$ be equivalent to $\mhA$, given as in Definition \ref{def:deng}, and $\mathcal{C}, \mathcal{D}$ and the unsmoothed prolongation $P$
  be as in Definition \ref{def:pwp}.
  Finally, for $I, J\in\mathcal{I}_c$, write $Q^{IJ}$ for $\sqrt{Q^{v_I\rightarrow v_J}}$ and $Q^{IJ\rightarrow ij} := Q^{ij}Q^{v_I\rightarrow v_i}Q^{JI} $. Then, with
  \[
    E^{IJ} := \sum_{v_i\in \mathcal{C}_I,v_j\in \mathcal{C}_J}Q^{IJ\rightarrow ij,T}E^{ij}Q^{IJ\rightarrow ij}
  \]
  and
  \[
    M^{I} := \sum_{v_i\in \mathcal{C}_I} Q^{v_I\rightarrow v_i, T}\left(M^i + \sum_{v_j\in \mcD} Q^{ij, T} E^{ij} Q^{ij}\right)Q^{v_I\rightarrow v_i},
  \]
  the auxiliary matrix $\mhA_c$, induced by these matrices as in Definition \ref{def:deng}, is equivalent to the coarse grid matrix,
  \[
    \norm{u}_{\mP^T\mA\mP} \sim \norm{u}_{\mhA_c},
  \]
  with constants that are no worse then those in the equivalence between $\mA$ and $\mhA$.
\end{theorem}
\begin{proof}
  The statement immediately follows from $\mP^T\mhA\mP = \mhA_c$. To see that, first consider the vertex contributions,
  \begin{align*}
    \sum_{v_i\in\mcV}\norm{(\mP\vof{u})_{i}}^2_{M^i} =& \sum_{C_I\in\mcC}\sum_{v_i\in C_I}\norm{(\mP\vof{u})_{i}}^2_{M^i} + \sum_{v_i\in\mcD}\norm{(\mP\vof{v})_i}^2_{M^i} \\
    =& \sum_{C_I\in\mcC}\sum_{v_i\in C_I}\norm{Q^{v_I\rightarrow v_i}\vof{u}_{I}}^2 + \sum_{v_i\in\mcD}\norm{0}^2\\
    =& \sum_{C_I\in\mcC}\vof{u}_I^T \left(\sum_{v_i\in C_I}Q^{v_I\rightarrow v_i, T}M^iQ^{v_I\rightarrow v_i}\right) \vof{u}_I.
  \end{align*}
  For the edge contributions, first see that 
  \[
    Q^{IJ\rightarrow ij}Q^{IJ} = Q^{ij}Q^{I\rightarrow i},
  \]
  which follows from elementary computations. With this there holds
  \begin{align*}
    \sum_{v_i\neq v_j\in\mcV} \norm{Q^{ij}(\mP\vof{u})_i-Q^{ji}(\mP\vof{u})_j}^2 = &
    \sum_{C_I\neq C_J\in\mcC}\sum_{v_i\in C_I} \sum_{v_j\in C_J} \norm{Q^{ij}(\mP\vof{u})_i-Q^{ji}(\mP\vof{u})_j}^2 +
    \sum_{C_I\in\mcC}\sum_{v_i\in C_I}\sum_{v_j\in\mcD}\norm{Q^{ij}(\mP\vof{u})_i-Q^{ji}(\mP\vof{u})_j}^2\\
    = &
    \sum_{C_I\neq C_J\in\mcC}\sum_{v_i\in C_I} \sum_{v_j\in C_J} \norm{Q^{ij}Q^{v_I\rightarrow v_i}\vof{u}_I-Q^{ji}Q^{v_J\rightarrow v_j}\vof{u}_J}^2 +
    \sum_{C_I\in\mcC}\sum_{v_i\in C_I}\sum_{v_j\in\mcD}\norm{Q^{ij}Q^{v_I\rightarrow v_i}\vof{u}_I}^2\\
    = &
    \sum_{C_I\neq C_J\in\mcC}\sum_{v_i\in C_I, v_j\in C_J} \norm{Q^{IJ\rightarrow ij}(Q^{IJ}\vof{u}_I-Q^{JI}\vof{u}_J)}^2 +
    \sum_{C_I\in\mcC}\sum_{v_i\in C_I}\sum_{v_j\in\mcD}\norm{\vof{u}_I}^2_{Q^{v_I\rightarrow v_i, T}Q^{ij, T}E^{ij}Q^{ij}Q^{v_I\rightarrow v_i}}\\
    = &
    \sum_{C_I\neq C_J\in\mcC}\sum_{v_i\in C_I, v_j\in C_J} \norm{Q^{IJ}\vof{u}_I-Q^{JI}\vof{u}_J}^2_{Q^{IJ\rightarrow ij, T}E^{ij}Q^{IJ\rightarrow ij}} +
    \sum_{C_I\in\mcC}\sum_{v_i\in C_I}\sum_{v_j\in\mcD}\norm{\vof{u}_I}^2_{Q^{v_I\rightarrow v_i, T}Q^{ij, T}E^{ij}Q^{ij}Q^{v_I\rightarrow v_i}}.
  \end{align*}
  Summing both sides of the above equalities, and defining $E^{IJ}$ and $M^I$ as proposed in the theorem gives
  \[
    \norm{\mP\vof{u}}_{\mhA}^2 = \norm{\vof{u}}_{\mhA_c}^2 \quad\forall\vof{u}\in(\mathbb{R}^k)^{|\mathcal{C}|}
  \]
  and thus
  \[
    \mA_c = \mP^T\mat{A}\mP \sim \mP^T\mat{\mhA}\mP = \mhA_c.
  \]
\end{proof}

In the following, we assume appropriate $E^{ij}$ and $M^i$ to be given on the finest level, and coarse edge- and vertex-matrices to be computed according to Theorem \ref{thm:cemats}.
As mentioned before, it seems that in practice, for many cases, very crude choices, based on the absolute values of off-diagonal entries,
suffice for the discretizations considered here.


%% file: crs2.tex
\newcommand\Qij{Q^{ij}}
\newcommand\QijT{Q^{ij, T}}
\newcommand\Qji{Q^{ji}}
\newcommand\QjiT{Q^{ji, T}}

\newcommand\Saijf{Q^{ji, T}\mhD_{ii}Q^{ji}}
\newcommand\Sajif{Q^{ij, T}\mhD_{jj}Q^{ij}}
\newcommand\Saij{\mtD^{ij}}
\newcommand\piSaij{\mtD^{ij, \dagger}}
\newcommand\Saji{\mtD^{ji}}
\newcommand\piSaji{\mtD^{ji, \dagger}}

\newcommand\Eilf{Q^{li, T}E^{il}Q^{li}}
\newcommand\Eil{\tilde{E}^{il}}
\newcommand\Ejlf{Q^{lj, T}E^{jl}Q^{lj}}
\newcommand\Ejl{\tilde{E}^{jl}}

\newcommand{\E}[1]{\bm{#1}}

The next goal is to introduce a coarsening algorithm that is directly based on the auxiliary matrix $\mhA$.
As the finite element matrix $\mA$ itself is not needed for this step, in the multi-level case, the entire hierarchy of auxiliary matrices, represented
by their topology of vertices and edges as well as edge and vertex matrices, can be built according to Theorem \ref{thm:cemats} before computing any coarse level matrices.
We therefore suggest the name \emph{Auxiliary Topology AMG} for the method.

First, we proof that a coarsening that produces a good two-grid method for $\mhA$ also produces one for $\mA$.
Lemma \ref{thm:ktg_generic} is a slight extension of a result from \cite{Falgout2004, mlfvass} that bounds the condition number of the two grid method by
$$
K_{TG} = \sup_{\vv} \frac{\vv^T\mM(\mI-\mtP(\mtP^T\mM\mtP)^{-1}\mtP^T\mM)\vv}{\vv^T\mA\vv}
$$
if a convergent and SPD smoother $\mM$ is used.
If $\mM$ is not symmetric, but equivalent to a symmetric matrix, it can be replaced in the above expression at the cost of an additional constant.
For a thorough discussion of this, see \cite{mlfvass}.
We take extra care to also cover the case where $\mhA$ is singular, as can happen for the elasticity problem.
\begin{lemma}\label{thm:ktg_generic}
  Let $\mhA$ be an auxiliary matrix as in Definition \ref{def:deng} and let the tentative prolongation $\mP$ be of the form of Definition \ref{def:pwp}.
  Let the columns of the matrix $\mE$ consist of an orthonormal basis of $\ker{\mhA}^{\perp}$ and the columns of $\mE_c$ consist of an orthonormal basis of $\ker{(\mP^T\mhA\mP)}^\perp$ where both have a block-structure as in Lemma \ref{thm:eet_is_diag}..
  Let the finite element matrix $\mA$ be equivalent to $\mtA = \mE^T\mhA \mE$,
  \begin{align}
    c_1 \mtA \leq \mA \leq c_2 \mtA \label{eq:matequiv}
  \end{align}
  for some $c_1, c_2>0$ and let $\mM$ be a convergent smoother for $\mA$ such that
  $$
  d_1 \widetilde{\mM} \leq \mD \leq d_2 \widetilde{\mM}
  $$
  for some $d_1, d_2>0$, where $\widetilde{\mM} = \mM(\mM+\mM^T-\mA)^{-1}\mM^T>0$ and $\mD$ is the diagonal of $\mA$.
  Let $\mtD$ be the diagonal of $\mtA$ and $\mhD = \mE\mtD\mE^T$ be the corresponding diagonal matrix of $\mhA$.
  Then the two-grid method defined by the prolongation $\mtP=\mE^T\mP\mE_c$ and the smoother $\mM$ produce an SPD coarse grid matrix $\mtP^T\mA\mtP$ and
  \begin{align}
    \mA \leq \mB_{TG} \leq \frac{c_2 d_2}{c_1} K_{TG} \mA, \label{eq:tgc_generic}
  \end{align}
  where
  \begin{align}
  K_{TG} = \sup_{\vv\in\ker{\mhA}^{\perp}}{\frac{\vv^T \mhD(\mI-\mP(\mP^T\mhD\mP)^{\dagger}\mP^T\mhD)\vv}{\vv^T\mhA\vv}}. \label{eq:ktg_generic}
  \end{align}
\end{lemma}
\begin{proof}
  Firstly, $\mtP^T\mA\mtP$ is positive definite due to the definition of $\mE_c$ and
  under these prerequisites, Theorem 3.19 and the following Corollary 3.20 from \cite{mlfvass} give an estimate
  $$
  \mA \leq \mB_{TG} \leq \tilde{K}_{TG} \mA
  $$
  where
  $$
  \tilde{K}_{TG} = d_2~\sup_{\vv} \frac{\vv^T\mD(\mI-\mtP(\mtP^T\mD\mtP)^{-1}\mtP^T\mD)\vv}{\vv^T\mA\vv}.
  $$
  It remains to show that
  $$
  \tilde{K}_{TG} \leq \frac{c_2}{c_1} K_{TG}.
  $$
  Because $\mA$ and $\mtA$ are equivalent according to (\ref{eq:matequiv}), the same equivalence holds between their diagonal parts,
  \begin{align*}
    c_1 \mtD \leq \mD \leq c_2 \mtD.
  \end{align*}
  Also observe that, due to the block-diagonal structure of $\mE$ there holds $\mE\mtD\mE^T=\mhD$, and of course $\mE^T\mE=\mI$.
  Now
  \begin{align*}
    \vv^T\mD(\mI-\mtP(\mtP^T\mD\mtP)^{-1}\mtP^T\mD)\vv &
     = \inf_{\vw} \norm{\vv - \mtP\vw}^2_{\mD} 
     \leq c_2~ \inf_{\vw} \norm{\vv - \mtP\vw}^2_{\mtD}, \\ 
      &= c_2~\inf_{\vw} \norm{\mE\vv - \mE\mtP\vw}^2_{\mhD}
      = c_2~\inf_{\vw} \norm{\mE\mE^T(\mE\vv - \mP\mE_c\vw)}^2_{\mhD}.
  \end{align*}
  The range of $\mI - \mE\mE^T$ is the kernel of $\mhA$ and Lemma \ref{thm:eet_is_diag} tells us that the matrix is diagonal,
  so its range is also in the kernel of $\mhD$, $(\mI - \mE\mE^T)\mhD = 0$. 
  Therefore, finally,
  \begin{align*}
    \vv^T\mD(\mI-\mtP(\mtP^T\mD\mtP)^{-1}\mtP^T\mD)\vv &
      \leq c_2~ \inf_{\vw} \norm{\mE\mE^T(\mE\vv - \mP\mE_c\vw)}^2_{\mhD}
      = c_2~\inf_{\vw} \norm{\mE\vv - \mP\mE_c\vw}^2_{\mhD} \\
      & = c_2~\inf_{\vw} \norm{\mE\vv - \mP\vw}^2_{\mhD}
      = c_2~(\mE\vv)^T \mhD(\mI-\mP(\mP^T\mhD\mP)^{\dagger}\mP^T\mhD)(\mE\vv),
  \end{align*}
  which is the wanted estimate for the enumerator in (\ref{eq:ktg_generic}).
  As for the denominator,
  \begin{align*}
    \vv^T\mA\vv\geq c_1 \vv^T\mtA\vv = c_1 (\mE\vv)^T\mhA(\mE\vv)
  \end{align*}
  shows, together with the above estimate that
  \begin{align*}
  \tilde{K}_{TG} \leq \frac{c_2}{c_1} \sup_{\vv} \frac{(\mE\vv)^T \mhD(\mI-\mP(\mP^T\mhD\mP)^{\dagger}\mP^T\mhD)(\mE\vv)}{(\mE\vv)^T\mhA(\mE\vv)} = \frac{c_2}{c_1} \sup_{\vv\in\ker{\mhA}^{\perp}}{\frac{\vv^T \mhD(\mI-\mP(\mP^T\mhD\mP)^{\dagger}\mP^T\mhD)\vv}{\vv^T\mhA\vv}}.
  \end{align*}
\end{proof}
We now proceed to show an estimate for $K_{TG}$ which should be seen as a generalization one originally formulated for scalar problems in \cite{Napov2010}, Theorem 3.2.
From this we will derive a criterion for agglomerates that implies boundedness of $K_{TG}$.
As checking this criterion can be computationally expensive in practice, in Section \ref{sec:ccm} we give simplified coarsening criteria for "matchings", that is, coarsenings where every agglomerate consists of only one or two vertices.
Based on these, in Section \ref{sec:spw}, a coarsening algorithm based on successive rounds of matching, first introduced in \cite{Napov2012}, will be modified
to work for both Problem \ref{prob:gg} and Problem \ref{prob:epseps}.
We stress that two-grid convergence is guaranteed even for problems that have deliberately been constructed to defeat coarsening algorithms based on traditional, scalar criteria, as demonstrated in Section \ref{sec:res}.

\newcommand{\mud}{\mu_{\mcD}}
\newcommand{\mua}{\mu_{g}^{\delta}}
\newcommand{\mub}{\mu_{p}^{\delta}}
\newcommand{\muc}{\mu_{s}}

\begin{theorem}\label{thm:mu1}
  For every $C\in\mcC$, let $\mhD_C$ be the diagonal block of $\mhD$ associated to vertices in $C$, and let $\mP_C\in(\mathbb{R}^{\kdim\times\kdim})^{|C|\times 1}$ be
  the the block of $\mP$ associated to $C$. Now, for $\delta\in\{0,1\}$, let $\mua(C)$ be the smallest constant such that for all $\vv_C\in(\mathbb{R}^{\kdim})^{|C|}$
  \begin{align}
  \inf_{\vw\in\mathbb{R}^{\kdim}} \norm{\vv_C - \mP_C\vw}_{\mhD_C}^2 \leq
  \mua(C)\left(\sum_{v_i\in C}\norm{\vv_i}_{M^i}^2 +
  \sum_{v_i\neq v_j\in C} \norm{\Qij\vv_i-\Qji\vv_j}^2_{E^{ij}} +
    \frac{\delta}{2}\sum_{v_j\notin C}\inf_{\vw\in\mathbb{R}^{\kdim}}\sum_{v_i\in C}\norm{\Qij\vv_i - \Qji\vw}^2_{E^{ij}}\right) \label{eq:mu_a}
  \end{align}
  and let
  \begin{align}
  \mu_{\mcD}(v_j) := \frac{\vv_i^T\mhD_{ii}\vv_i}{\vv_i^T M^{i} \vv_i}. \label{eq:mu_d}
  \end{align}
  Then there holds
  \begin{align}
    \KTG \leq \max\{\sup_{v_i\in\mcD}\mu_{\mcD}(v_i), \sup_{C\in\mcC} \mua(C)\}. \label{eq:mu_k}
  \end{align}
\end{theorem}
\begin{proof}
  First, we find an estimate from below for the enumerator in (\ref{eq:ktg_generic}) by omitting edge-terms coming from connections between
  vertices in $\mcD$, scaling those coming from connections between
  vertices in $\mcD$ and vertices in some agglomerate $\mcC$ by one half and reordering the rest of the terms.
  \begin{align*}
    \norm{\vv}_{\mhA}^2 & = \sum_{v_i\in\mcV} \norm{\vv_i}_{M^i}^2 + \sum_{v_i\neq v_j\in\mcV}\norm{\Qij\vv_i-\Qji\vv_j}^2_{E^{ij}} \\
                        & \geq \sum_{v_i\in\mcD} \norm{\vv_i}_{M^i}^2 + \sum_{C\in\mcC}\left( \sum_{v_i\in C} \norm{\vv_i}_{M^i}^2 +
                          \sum_{v_i\neq v_j\in C}\norm{\Qij\vv_i-\Qji\vv_j}^2_{E^{ij}} +
                          \frac{1}{2}\sum_{v_i\in C}\sum_{v_j\notin C}\norm{\Qij\vv_i-\Qji\vv_j}^2_{E^{ij}}\right) \\
                        & \geq \sum_{v_i\in\mcD} \norm{\vv_i}_{M^i}^2 + \sum_{C\in\mcC}\left( \sum_{v_i\in C} \norm{\vv_i}_{M^i}^2 +
                          \sum_{v_i\neq v_j\in C}\norm{\Qij\vv_i-\Qji\vv_j}^2_{E^{ij}} +
                          \frac{1}{2}\sum_{v_j\notin C}\inf_{\vw}\sum_{v_i\in C}\norm{\Qij\vv_i-\Qji\vw}^2_{E^{ij}}\right) \\
                        & \geq \sum_{v_i\in\mcD} \norm{\vv_i}_{M^i}^2 + \sum_{C\in\mcC}\left( \sum_{v_i\in C} \norm{\vv_i}_{M^i}^2 +
                          \sum_{v_i\neq v_j\in C}\norm{\Qij\vv_i-\Qji\vv_j}^2_{E^{ij}} +
                          \frac{\delta}{2}\sum_{v_j\notin C}\inf_{\vw}\sum_{v_i\in C}\norm{\Qij\vv_i-\Qji\vw}^2_{E^{ij}}\right).
  \end{align*}
  On the other hand
  \begin{align*}
    \inf_{\vw}\norm{\vv-\mP\vw}_{\mhD}^2 & = \inf_{\vw}\sum_{i}\norm{\vv_i - (\mP\vw)_i}^2_{\mhD_{ii}} =
    \sum_{v_i\in\mcD}\norm{\vv_i}_{\mhD_{ii}}^2 + \inf_{\vw} \sum_{C\in\mcC}\sum_{v_i\in C}\norm{\vv_i - \mP_{iJ}\vw_J}^2_{\mhD_{ii}} \\
    & = \sum_{v_i\in\mcD}\norm{\vv_i}_{\mhD_{ii}}^2 + \sum_{C\in\mcC} \inf_{\vw}  \sum_{v_i\in C}\norm{\vv_i - \mP_{iJ}\vw}^2_{\mhD_{ii}} \\
    & = \sum_{v_i\in\mcD}\norm{\vv_i}_{\mhD_{ii}}^2 + \sum_{C\in\mcC} \inf_{\vw}  \norm{\vv_C - \mP_C\vw}_{\mhD_C}^2,
  \end{align*}
  and thus seperating the terms in the above sums proves (\ref{eq:mu_k}).
\end{proof}
Computing $\mua(C)$ involves finding the largest eigenvalue of the generalized eigenvalue problem
$$
\mhM_C\vv = \lambda \mhA_C\vv,
$$
where $\mhM_C, \mhA_C\in(\mathbb{R}^{\kdim\times\kdim})^{|C|\times |C|}$ are the matrices that induce the semi-norms on the left and right hand side, respectively.
The role of $\delta$ in (\ref{eq:mu_a}) is important for elasticity, where $\delta=0$ leads to a criterion that is too restrictive in practice.
The reason is that $\mhA_C$ then often has a lower rank than $\mhM_C$ due to the fact that both edge and vertex matrices can be singular, and indeed vertex matrices are often zero, and the inequality can not hold for any finite $\mu_{g}^{0}(C)$.
This results in a coarsening algorithm producing very large coarse grids and eventually, when applied on multiple levels, completely stalling. 
Setting $\delta=1$ fixes the rank deficiency in enough cases to obtain good coarsening rates.
\begin{remark}
  When $\delta=1$ in (\ref{eq:mu_a}), the local energy on the right hand side is induced by a (generalized, because edge matrices can be singular)
  Schur complement matrix.
  This is similar to the way strength of connection between two vertices is computed in AMGm, \cite{CMIII}, based on a generalized Schur complement of a
  so called computational molecule, a small sub-assembled diagonal block of the auxiliary matrix, assembled from edge contributions between the two vertices and their common
  neighbors.
\end{remark}
As already proposed in \cite{Napov2012} for their coarsening algorithm for scalar problems, alternatively, in order to simply check that
$\mua(C)\leq \mu$ for some fixed parameter $\mu$, the matrix $\mu \mhA_C - \mhM_C$ can be checked for positive semi-definiteness via computing a Cholesky factorization.
It should be noted that when $\delta=1$, because no connections between outside neighbors are included, the right hand side matrix
can be computed efficiently by eliminating neighbors separately instead of all at once.

\subsection{Coarsening criteria for matchings}\label{sec:ccm}
Lemma \ref{thm:mu2} introduces a measure that bounds $\mua$ from above. While cheaper to compute, it is only applicable to pairs of vertices.
It reduces the eigenvalue problem that has to be solved from size $2\kdim$ to size $\kdim$ which makes it completely feasible for an eventual coarsening algorithm to exactly compute it, instead of just checking whether it is bounded by some constant.
\begin{lemma}\label{thm:mu2}
  Let $C=\{v_i, v_j\}$ be an agglomerate of size 2 and let $\mcN$ be the set of
  neighboring vertices common to $v_i$ and $v_j$. 
  For $\delta\in\{0,1\}$, let $\mub(v_i, v_j)$ be the largest eigenvalue of the generalized eigenvalue problem
  \begin{align}\label{eq:mub_ineq}
  H(\QjiT\mhD_{ii}\Qji, \QijT\mhD_{jj}\Qij) = \lambda \left( E^{ij} + \frac{\delta}{2}\sum_{l\in\mcN}\QjiT Q^{v_i\rightarrow v_l, T}H(Q^{li, T}E^{il}Q^{li}, Q^{lj, T}E^{jl}Q^{lj})Q^{v_i\rightarrow v_l}\Qji \right)
  \end{align}
  where $H(A, B) = A(A+B)^{\dagger}B$ is the harmonic mean.
  Then there holds
  $$
  \mua(\{v_i, v_j\}) \leq \mub(v_i, v_j).
  $$
\end{lemma}
\begin{proof}
  Let $v_H$ be the coarse vertex representing the agglomerate $C$ (for which there holds $(\mP_{C}\vw)_{l}=Q^{v_H\rightarrow v_l}\vw$).
  Set $Q^{Hh} := \Qij Q^{v_H\rightarrow v_i}$ and observe that $ \Qij Q^{v_H\rightarrow v_i} = \Qji Q^{v_H\rightarrow v_j}$.
  With this we show the following identity for the left hand side of (\ref{eq:mu_a}):
  \begin{align*}
    \inf_{\vw\in\mathbb{R}^{\kdim}} \norm{\vv_C-\mP_C\vw}_{\mhD_C}^2 & =
    \inf_{\vw\in\mathbb{R}^{\kdim}} \norm{\vv_i-Q^{v_H\rightarrow v_i}\vw}_{\mhD_{ii}}^2
  + \norm{\vv_j-Q^{v_H\rightarrow v_j}\vw}_{\mhD_{jj}}^2 \\
    & = \inf_{\vw\in\mathbb{R}^{\kdim}} \norm{\Qji\Qij(\vv_i-Q^{v_H\rightarrow v_i}\vw)}_{\mhD_{ii}}^2
  + \norm{\Qij\Qji\left(\vv_j-\vw\right)}_{\mhD_{jj}}^2 \\
    & = \inf_{\vw\in\mathbb{R}^{\kdim}} \norm{\Qji\left(\Qij\vv_i-\Qij Q^{v_H\rightarrow v_i}\vw\right)}_{\mhD_{ii}}^2
  + \norm{\Qij\left(\Qji\vv_j-\Qji Q^{v_H\rightarrow v_j}\vw\right)}_{\mhD_{jj}}^2 \\
    & = \inf_{\vw\in\mathbb{R}^{\kdim}} \norm{\Qij\vv_i-Q^{Hh}\vw}_{\QjiT\mhD_{ii}\Qji}^2
  + \norm{\Qji\vv_j-Q^{Hh}\vw}_{\QijT\mhD_{jj}\Qij}^2 \\
    & = \inf_{\vw\in\mathbb{R}^{\kdim}} \norm{\Qij\vv_i-\vw}_{\QjiT\mhD_{ii}\Qji}^2
  + \norm{\Qji\vv_j-\vw}_{\QijT\mhD_{jj}\Qij}^2.
  \end{align*}
  For readability, set $\Saij=\Saijf$ and $\Saji=\Sajif$. Elementary calculations show that the infimum is obtained for
  $$\vw =  \left(\Saij + \Saji\right)^{\dagger}\left(\Saij \Qij\vv_i + \Saji \Qji\vv_j\right).$$
  For a short detour, consider two symmetric, nonnegative matrices $A$ and $B$. Clearly, $\ker(A+B)=\ker A\cap\ker B$, so
  $\ker(A)^\perp \subseteq \ker(A+B)^\perp$,
  and thus there holds
  \begin{align}
  \norm{x}_{A} = \norm{(A+B)^{\dagger}(A+B)x}_{A}. \label{eq:sum_dagger_sum}
  \end{align}
  Now, inserting the solution for $\vw$ in the terms in the infimum above gives
  \begin{align*}
    \norm{\Qij\vv_i-\vw}_{\Saij}^2 & = \norm{\Qij\vv_i-\left(\Saij + \Saji\right)^{\dagger}\left(\Saij \Qij\vv_i + \Saji \Qji\vv_j\right)}_{\Saij}^2 \\
    & = \norm{\left(\Saij + \Saji\right)^{\dagger}\left(\Saij + \Saji\right)\left(\Qij\vv_i-\left(\Saij + \Saji\right)^{\dagger}\left(\Saij \Qij\vv_i + \Saji \Qji\vv_j\right)\right)}_{\Saij}^2 \\
    & = \norm{\left(\Saij + \Saji\right)^{\dagger}\left(\left(\Saij + \Saji\right)\Qij\vv_i - \Saij \Qij\vv_i - \Saji \Qji\vv_j\right)}_{\Saij}^2 \\
    & = \norm{\left(\Saij + \Saji\right)^{\dagger}\Saji\left(\Qij\vv_i -  \Qji\vv_j\right)}_{\Saij}^2 \\
    & = \norm{\Qij\vv_i -  \Qji\vv_j}_{\Saji\left(\Saij + \Saji\right)^{\dagger}\Saij\left(\Saij + \Saji\right)^{\dagger}\Saji}^2 \intertext{and, analogously,}
    \norm{\Qji\vv_i-\vw}_{\Saji}^2 & = \ldots =
    \norm{\Qij\vv_i -  \Qji\vv_j}_{\Saij\left(\Saij + \Saji\right)^{\dagger}\Saji\left(\Saij + \Saji\right)^{\dagger}\Saij}^2.
  \end{align*}
  As further elementary calculations show, for any two symmetric, nonnegative matrices $A, B$
  \begin{align}
    A(A+B)^{\dagger}B(A+B)^{\dagger}A + B(A+B)^{\dagger}A(A+B)^{\dagger}B &= A(A+B)^{\dagger}B \label{eq:harm_pinv}
  \end{align}
  and thus
  \begin{align*}
    \inf_{\vw\in\mathbb{R}^{\kdim}} \norm{\vv_C-\mP_C\vw}_{\mhD_C}^2 & = \norm{\Qij\vv_i -  \Qji\vv_j}_{\Saij\left(\Saij + \Saji\right)^{\dagger}\Saji}^2.
  \end{align*}
  As for the right hand side of (\ref{eq:mu_a}), we need to identify
  $$
  \sum_{v_l\notin C}\inf_{\vw\in\mathbb{R}^{\kdim}}\sum_{v_k\in C}\norm{Q^{kl}\vv_k - Q^{lk}\vw}^2_{E^{kl}}.
  $$
  All terms for vertices $v_l\notin\mcN$ (that are neighbor of only one vertex in $V$) are zero, which is attained for
  the choice $\vw = Q^{v_k\rightarrow v_l}\vv_k$. It is therefore enough to sum over all $v_l\in\mcN$:
  $$
  \sum_{v_l\notin C}\inf_{\vw\in\mathbb{R}^{\kdim}}\sum_{v_k\in C}\norm{Q^{kl}\vv_k - Q^{lk}\vw}^2_{E^{kl}} =
  \sum_{v_l\in \mcN}\inf_{\vw\in\mathbb{R}^{\kdim}}\norm{Q^{il}\vv_i - Q^{li}\vw}^2_{E^{il}} + \norm{Q^{jl}\vv_j - Q^{lj}\vw}^2_{E^{jl}}.
  $$
  Similarly to before, with $\Eil=\Eilf$ and $\Ejl=\Ejlf$, elementary calculations show that the infimum
  \begin{align*}
    \inf_{\vw\in\mathbb{R}^{\kdim}}\norm{Q^{il}\vv_i - Q^{li}\vw}^2_{E^{il}} + \norm{Q^{jl}\vv_j - Q^{lj}\vw}^2_{E^{jl}}
  \end{align*}
  is attained for
  $$
  \vw = \left(\Eil + \Ejl\right)^{\dagger}\left( Q^{li,T}E^{li}Q^{il}\vv_i + Q^{lj,T}E^{lj}Q^{jl}\vv_j\right).
  $$
  Inserting this into one of the terms in the infimum above gives
  \begin{align*}
    \norm{Q^{il}\vv_i - Q^{li}\vw}^2_{E^{il}} & =
\norm{Q^{il}\vv_i - Q^{li}\left(\Eil + \Ejl\right)^{\dagger}\left( Q^{li,T}E^{li}Q^{il}\vv_i + Q^{lj,T}E^{lj}Q^{jl}\vv_j\right)}^2_{E^{il}} \\
    & = \norm{Q^{li}\left(Q^{v_i\rightarrow v_l}\vv_i - \left(\Eil + \Ejl\right)^{\dagger}\left( Q^{li,T}E^{li}Q^{il}\vv_i + Q^{lj,T}E^{lj}Q^{jl}\vv_j\right)\right)}^2_{E^{il}} \\
                                              & = \norm{Q^{v_i\rightarrow v_l}\vv_i - \left(\Eil + \Ejl\right)^{\dagger}\left( Q^{li,T}E^{li}Q^{il}\vv_i + Q^{lj,T}E^{lj}Q^{jl}\vv_j\right)}^2_{\Eil}. \\
    \intertext{Applying (\ref{eq:sum_dagger_sum}) allows for further simplification:}
    \norm{Q^{il}\vv_i - Q^{li}\vw}^2_{E^{il}} & = \norm{\left(\Eil + \Ejl\right)^{\dagger}\left(\left(\Eil + \Ejl\right)Q^{v_i\rightarrow v_l}\vv_i - \left( Q^{li,T}E^{li}Q^{il}\vv_i + Q^{lj,T}E^{lj}Q^{jl}\vv_j\right)\right)}^2_{\Eil} \\
    & = \norm{\left(\Eil + \Ejl\right)^{\dagger}\left(\Eil Q^{v_i\rightarrow v_l}\vv_i + \Ejl Q^{v_i\rightarrow v_l}\vv_i - Q^{li,T}E^{li}Q^{il}\vv_i - Q^{lj,T}E^{lj}Q^{jl}\vv_j\right)}^2_{\Eil} \\
    & = \norm{\Eil Q^{v_i\rightarrow v_l}\vv_i + \Ejl Q^{v_i\rightarrow v_l}\vv_i - Q^{li,T}E^{li}Q^{li}Q^{il}Q^{il}\vv_i - Q^{lj,T}E^{lj}Q^{jl}\vv_j}^2_{\left(\Eil + \Ejl\right)^{\dagger}\Eil\left(\Eil + \Ejl\right)^{\dagger}} \\
    & = \norm{\Ejl Q^{v_i\rightarrow v_l}\vv_i - Q^{lj,T}E^{lj}Q^{lj}Q^{jl}Q^{jl}\vv_j}^2_{\left(\Eil + \Ejl\right)^{\dagger}\Eil\left(\Eil + \Ejl\right)^{\dagger}} \\
    & = \norm{\Ejl\left(Q^{v_i\rightarrow v_l}\vv_i - Q^{v_j\rightarrow v_l}\vv_j\right)}^2_{\left(\Eil + \Ejl\right)^{\dagger}\Eil\left(\Eil + \Ejl\right)^{\dagger}} \\
    & = \norm{Q^{v_i\rightarrow v_l}\Qji\Qij\vv_i - Q^{v_j\rightarrow v_l}\Qij\Qji\vv_j}^2_{\Ejl\left(\Eil + \Ejl\right)^{\dagger}\Eil\left(\Eil + \Ejl\right)^{\dagger}\Ejl} \\
    & = \norm{Q^{v_i\rightarrow v_l}\Qji\left(\Qij\vv_i - \Qji\vv_j\right)}^2_{\Ejl\left(\Eil + \Ejl\right)^{\dagger}\Eil\left(\Eil + \Ejl\right)^{\dagger}\Ejl} \\
    & = \norm{\Qij\vv_i - \Qji\vv_j}^2_{\QjiT Q^{v_i\rightarrow v_l, T}\Ejl\left(\Eil + \Ejl\right)^{\dagger}\Eil\left(\Eil + \Ejl\right)^{\dagger}\Ejl Q^{v_i\rightarrow v_l}\Qji}. \\
  \intertext{Analogous computations for the second term give}
    \norm{Q^{jl}\vv_i - Q^{li}\vw}^2_{E^{jl}}
    & = \norm{\Qij\vv_i - \Qji\vv_j}^2_{\QjiT Q^{v_i\rightarrow v_l, T}\Eil\left(\Ejl + \Eil\right)^{\dagger}\Ejl\left(\Ejl + \Eil\right)^{\dagger}\Eil Q^{v_i\rightarrow v_l}\Qji}.
  \end{align*}
  Here, we used
  $$
  Q^{v_i\rightarrow v_l}\Qji = Q^{v_j\rightarrow v_l}\Qij,
  $$
  which can easily be seen by
  $$
  Q^{v_i\rightarrow v_l}\Qji\left(Q^{v_j\rightarrow v_l}\Qij\right)^{-1} = Q^{v_i\rightarrow v_l}\Qji\Qji Q^{v_l\rightarrow v_j} = Q^{v_i\rightarrow v_l}Q^{v_j\rightarrow v_i}Q^{v_l\rightarrow v_j} = I.
  $$
  Summing these two and using (\ref{eq:harm_pinv}) we have found the infimum,
  \begin{align*}
    \inf_{\vw\in\mathbb{R}^{\kdim}}\norm{Q^{il}\vv_i - Q^{li}\vw}^2_{E^{il}} + \norm{Q^{jl}\vv_j - Q^{lj}\vw}^2_{E^{jl}}
    & = \norm{\Qij\vv_i - \Qji\vv_j}^2_{\QjiT Q^{v_i\rightarrow v_l, T}H(\Ejl, Q^{li, T}E^{ij}Q^{li})Q^{v_i\rightarrow v_l}\Qji}. \\
  \end{align*}
  Summing over all $v_l\in\mcN$,
  $$
    H(\QjiT\mhD_{ii}\Qji, \QijT\mhD_{jj}\Qij) \leq \mub(v_i, v_j) \left( E^{ij} + \frac{\delta}{2}\sum_{l\in\mcN}\QjiT Q^{v_i\rightarrow v_l, T}H(\Eilf, \Ejlf)Q^{v_i\rightarrow v_l}\Qji \right)
  $$
  is equivalent to
  \begin{align*}
  \inf_{\vw\in\mathbb{R}^{\kdim}} \norm{\vv_C - \mP_C\vw}_{\mhD_C}^2 &\leq \mub(v_i, v_j)\left(
  \sum_{v_i\neq v_j\in C} \norm{\Qij\vv_i-\Qji\vv_j}^2_{E^{ij}} +
\frac{\delta}{2}\sum_{v_l\notin C}\inf_{\vw\in\mathbb{R}^{\kdim}}\sum_{v_k\in C}\norm{Q^{kl}\vv_k - Q^{lk}\vw}^2_{E^{kl}}\right),
  \intertext{which implies}
    \inf_{\vw\in\mathbb{R}^{\kdim}} \norm{\vv_C - \mP_C\vw}_{\mhD_C}^2 
    &\leq \mub(v_i, v_j)\left( \sum_{v_i\in C} \norm{\vv_j}_{M^i}^2 + 
  \sum_{v_i\neq v_j\in C} \norm{\Qij\vv_i-\Qji\vv_j}^2_{E^{ij}} +
\frac{\delta}{2}\sum_{v_l\notin C}\inf_{\vw\in\mathbb{R}^{\kdim}}\sum_{v_k\in C}\norm{Q^{kl}\vv_k - Q^{lk}\vw}^2_{E^{kl}}\right).
  \end{align*}
  Finally, as $\mua(\{v_i, v_j\})$ is the smallest constant for which this last inequality holds,
  $$
  \mua(\{v_i, v_j\}) \leq \mub(v_i, v_j).
  $$
\end{proof}

\begin{figure}
  \centering
  \includegraphics[scale=0.5]{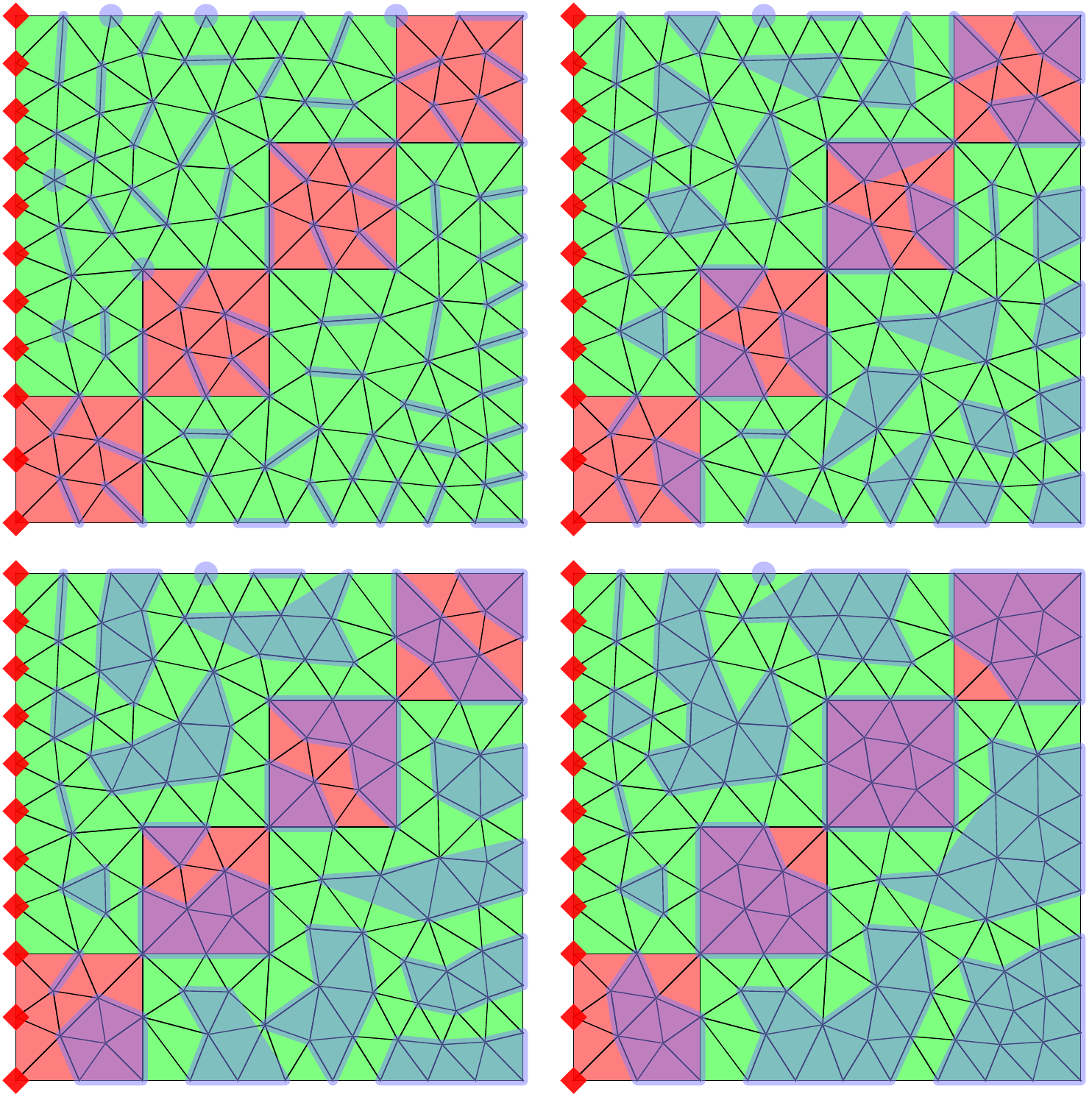}
  \caption{Agglomerates formed by Algorithm \ref{alg:pwagg} for the linearized elasticity problem, (\ref{eq:epseps}), on $\Omega=(0,1)^2$.
    In the green parts of the domain $\mu=\lambda=1$ and $\mu=\lambda=10^4$ in the red parts.
    Dirichlet conditions are prescribed on the left side.
    Note that the algorithm does not group up vertices in different ``boxes''.
    For the two-grid method with Gauss-Seidel smoother using this coarsening, $\kappa(\mB_{TG}^{-1}\mA)\approx 2.84$.
  }
  \label{fig:2dc_rob}
\end{figure}

\begin{figure}
  \centering
  \includegraphics[scale=0.5]{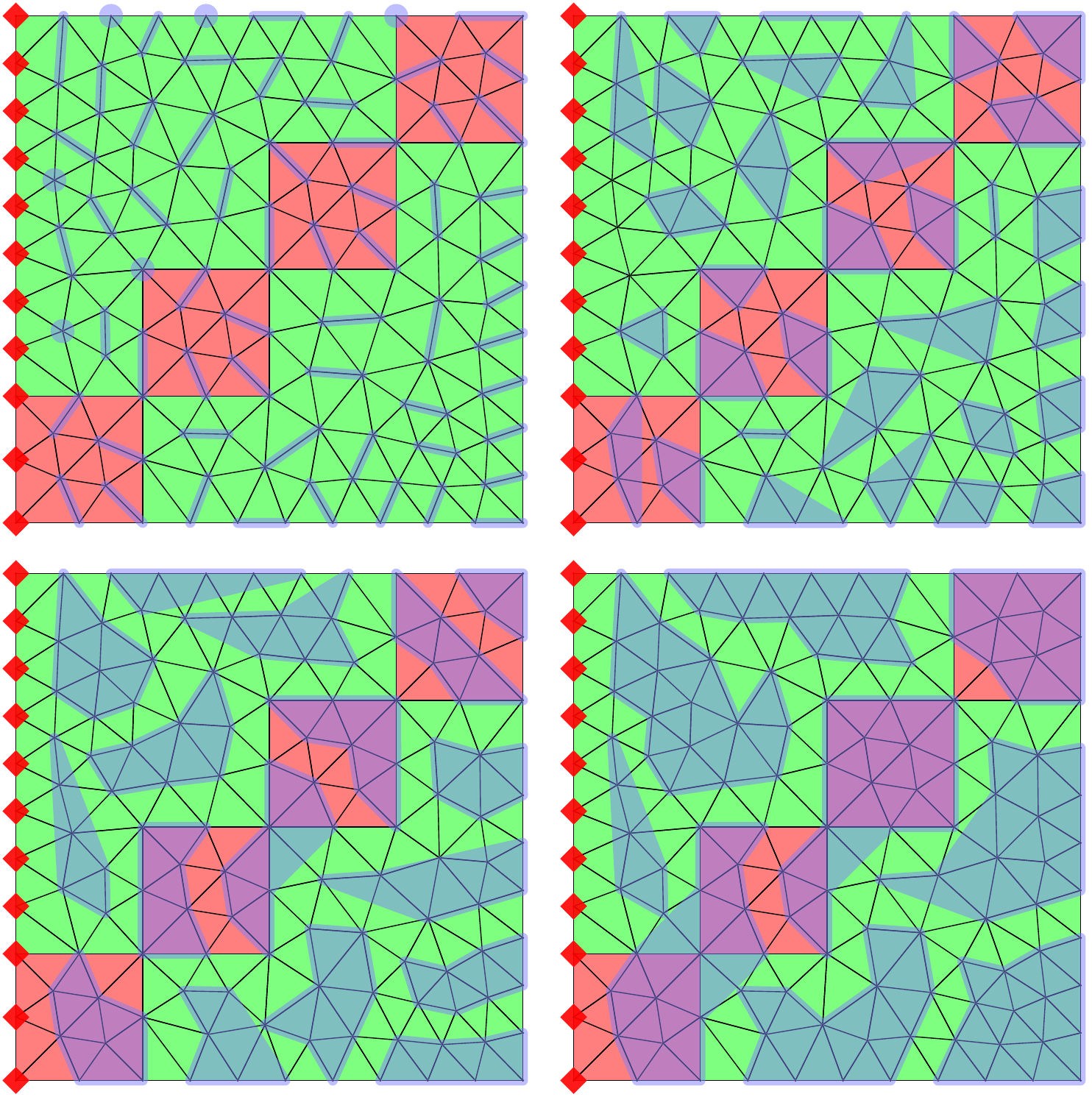}
  \caption{Agglomerates formed by Algorithm \ref{alg:pwagg} for the same problem as in Figure \ref{fig:2dc_rob}
    when only using the scalar measure $\muc$.
    That means $\mua$ and $\mub$ are not computed in lines \ref{algline:mub} and \ref{algline:mua} of the algorithm.
    Note that vertices from different ``boxes'' are grouped up already in the second step.
    This time $\kappa(\mB_{TG}^{-1}\mA)\approx 270$.
  }
  \label{fig:2dc_norob}
\end{figure}

\begin{remark}
  In fact, as can be seen in the last step of the proof, $\mub(v_i, v_j) = \mua(\{v_i, v_j\})$ if $M^i = M^j = 0$.
\end{remark}
In \cite{Napov2012}, a similar simplification for pairs of vertices is made for the scalar case which results in a scalar expression that is slightly different from
$\mub$ in that is also includes diagonal contributions analogous to $M^i, M^j$.
For elasticity problems, where $\kdim>1$, it would be convenient to find an additional simplification of $\mub$ which results in some kind of scalar criterion. 
During the coarsening process, it could be used to cheaply pick a candidate for agglomeration among viable neighbors.
As long as that choice is then confirmed by $\mub$, it does not even necessarily need to be an upper bound for $\mub$.
We use a heuristic measure inspired by a criterion that is used in classical scalar AMG methods for M-matrices to gauge the strength of connection
between two degrees of freedom $i$ and $j$, based on entries of the finite element matrix $\mA$ (see also, for example, the overview paper \cite{amg_ov_Xu2017}),
$$
s_{ij} = \frac{|\mA_{ij}|}{\min\{\max_{k\neq i}\{|\mA_{ik}|\}, \max_{k\neq j}\{|\mA_{jk}|\}\}}.
$$
Instead of using entries of $\mA$, or even $\mhA$, it is even simpler to use vertex- and edge-matrices $M^{i}, E^{ij}$, which are guaranteed to have nonnegative trace:
\begin{definition}\label{def:muc}
  With the trace of a matrix $\trace{A}:=\sum_iA_{ii}$, for two vertices $v_i,v_j\in\mcV$, the scalar measure $\muc(v_i, v_j)$ is defined as
  \begin{align*}
  \muc(v_i, v_j) := \frac{ \sqrt{\max\{\trace{M^{i}}, \max_{l\neq i}{\{\trace{E^{il}}\}\}}\cdot\max\{\trace{M^{j}}, \max_{l\neq j}{\{\trace{E^{jl}}}\}\}} }{\trace{E^{ij}}} \label{eq:muc}
  \end{align*}
\end{definition}
Note however, that $\muc$ is in no way \textit{guaranteed} to produce high quality agglomerates.
As an example, consider Figure \ref{fig:2dc_rob}, where the combination of $\mua, \mub,$ and $\muc$ leads to a well conditioned two-grid method, and Figure \ref{fig:2dc_norob},
where only $\muc$ is considered while $\mua$ and $\mub$ are ignored and thus a much worse condition number is attained.
It should be said, however, that this example was specifically designed to make proper coarsening difficult.
In many cases, however, $\muc$ suffices in practice.

\subsection{A successive matching algorithm}\label{sec:spw}

The coarsening algorithm presented here (Algorithm \ref{alg:pwagg}) is based on the successive pairwise aggregation algorithm from \cite{Napov2012} which matches up vertices in multiple passes.

In the first pass, vertices are matched up and coarse grid vertex and edge matrices are computed according to Theorem \ref{thm:cemats} in line \ref{algline:calccemats1}.
As no agglomerates that contain more than two vertices can be formed, this can be done based on $\mub, \muc$ alone.
On subsequent passes, the (coarse grid) vertices that are paired up represent agglomerates formed by the previous pass.
As now also bigger agglomerates exist, 
also $\mua$ (which still depends on fine grid matrices $M^i, E^{ij}$) has to be considered.
Note that on later passes, in line \ref{algline:calccemats2}, edge and vertex matrices for the next one are computed based on previous coarse grid matrices and not initial fine grid matrices.

In the matching process itself, candidates have to be selected among unmatched neighbors.
As Theorem \ref{thm:mu1} shows that the measure $\mua$ goes directly into the estimate for $K_{TG}$, it would be preferrable to pick the neighbor that minimizes $\mua$.
On the first pass, also the neighbor that minimizes $\mub$ (which is an upper bound for $\mua$) would a reasonable choice.
However, on later passes, computing $\mua$ for every single neighbor would be expensive, as increasingly large eigenvalue problems have to be solved.
Therefore, as suggested in \cite{Napov2012}, instead of the "optimal" neighbor, a candidate is chosen that maximizes either $\mub$ or $\muc$ (computed from coarse grid matrices $E^{IJ}, M^{I}$).
If the initial choice is made based on $\muc$, it is first confirmed by also computing $\mub$.
Where necessary (that is, if the resulting agglomerate would contain more than three vertices), this is then finally confirmed by checking that $\mua<\sigma$, which can be done, as mentioned before, by computing a Cholesky factorization.
The advantage of basing the initial choice on $\muc$ is that the more expensive $\mub$ ideally has to be computed for only one candidate.
On the other hand, we might get a better initial candidate from $\mub$, as it is an upper bound to $\mua$ on the firt pass and an upper bound to what would be $\mua$ computed from the coarse grid matrix on later passes.
However, in practice, we have not overseved that have a considerable effect.

One detail that has a surprisingly large effect on the quality of generated agglomerates is the choice of $C_i$ in line \ref{algline:pickci}.
As again suggested in \cite{Napov2012}, on the first pass, a Cuthill-McKee permutation of fine grid vertices is computed and vertices are iterated through in that manner.
On subsequent passes, vertices are iterated through in the reverse order they have been added to the coarse grid on the previous pass.
The same ordering is also used to break potential ties in line \ref{algline:pickmincj}.

{
  \newcommand{\Algphase}[1]{%
    \vspace*{-.7\baselineskip}\Statex\hspace*{\dimexpr-\algorithmicindent-2pt\relax}\rule{\textwidth}{0.4pt}%
    \Statex\hspace*{-\algorithmicindent}\textbf{#1}%
    \vspace*{-.7\baselineskip}\Statex\hspace*{\dimexpr-\algorithmicindent-2pt\relax}\rule{\textwidth}{0.4pt}%
  }
  \algrenewcommand\algorithmicrequire{\textbf{Input:}}
  \algrenewcommand\algorithmicensure{\textbf{Output:}}
  \algnotext{EndFor}
  \algnotext{EndProcedure}
  \algnotext{EndIf}
  \algnotext{EndWhile}
  \begin{algorithm}
    \caption{successive matching algorithm}
    \label{alg:pwagg}
    \begin{algorithmic}[1] 
      \Statex{}
      \Require {\text{ } \newline set of vertices $\mcV$ and Dirichlet vertices $\mcV_D$ \newline vertex- and edge-matrices $M^i, E^{ij}$, represented by $\mhA$ \newline threshold $\sigma$ \newline $\delta\in\{0,1\}$ determines the choice between $\mu_{g}^{0},\mu_{p}^{0}$ and $\mu_{g}^{1},\mu_{p}^{1}$ \newline $R\in\{0,1\}$ determines whehter $\mub$ or $\muc$ is used to pick a candidate neighbour \newline number of rounds $N$}
      \Ensure {\text{ }\newline set of agglomerates $\mcC$, $\mcD$}
      \vspace*{-.5\baselineskip}\Statex\hspace*{\dimexpr-\algorithmicindent-2pt\relax}\rule{\textwidth}{0.4pt}\vspace*{-\baselineskip}
      \Statex{}
      \State {Set $\mcD = \mcV_D\cup\{ v_i : \mud(v_i)<\sigma\}$}
      \State {Set $\mcC = \{\{v_i\} : v_i\in\mcC\setminus\mcD\}$}
      \For{$k=1,\ldots, N$}
      \State{$\mcT$ = \Call{match\_vertices}{$\sigma$, $\delta$, $R$, $\mcC$} }
      \If{$k<N$}
      \If{$k==0$}
      \State{Set up $E^{IJ}, M^{I}$ from $E^{ij}, M^i$ according to Theorem \ref{thm:cemats}, using agglomerates $\mcT$ and Dirichlet set $\mcD$ } \label{algline:calccemats1}
      \Else
      \State{Set up new $E^{IJ}, M^{I}$ from $E^{IJ}, M^{I}$ according to Theorem \ref{thm:cemats}, using agglomerates $\mcT$ } \label{algline:calccemats2}
      \EndIf
      \EndIf
      \State{Set $\mcC = \{\bigcup_{C\in T}C : T\in\mcT\}$}
      \EndFor
      \State{\Return $\mcC, \mcD$}
      \Statex{$ $}
      \Algphase{a single round of matching}
      \Procedure{match\_vertices}{$\sigma$, $\delta$, $R$, $\mcT$}
      \State {Set $\mcC = \emptyset$, $\mcD = \emptyset$, $\mcS=\emptyset$}
      \While{$\mcT\neq\emptyset$}
      \State{Select $C_i\in\mcT$} \label{algline:pickci}
      \State{Set $\mcT=\mcT\setminus\{C_i\}$}
      \State{Set $\mcN_{s} = \{C_j \in\mcT: \muc(C_i, C_j)<\sigma\}$} \label{algline:muc}
      \While{$\mcN_{s}\neq\emptyset$}
      \If{$R$}
      \State{Select $C_j\in\mcN_{s}$ with minimal $\muc(C_i, C_j)$} \label{algline:pickmincj}
      \Else
      \State{Select $C_j\in\mcN_{s}$ with minimal $\mub(C_i, C_j)$} \label{algline:pickmincj2}
      \EndIf
      \State{Set $\mu=\mub(C_i, C_j)$} \label{algline:mub}
      \If{$\mu<\sigma$ \textbf{and} $\left(|C_i\cup C_j|=2 \textbf{ or } \mua(C_i\cup C_j)<\sigma\right)$} \label{algline:mua}
      \State{Set $\mcC = \mcC \cup \{ \{C_i, C_j\}\}$}
      \State{Set $\mcT = \mcT \setminus \{C_j\}$}
      \State{\textbf{break}}
      \Else
      \State{Set $\mcN_{s} = \mcN_{s}\setminus\{C_j\}$}
      \EndIf
      \EndWhile
      \If{$\mcN_s = \emptyset$}
      \State{Set $\mcS = \mcS\cup\{\{C_i\}\}$}
      \EndIf
      \EndWhile
      \State{Set $\mcC = \mcC \cup \mcS$}
      \State{\Return{$\mcC$}}
      \EndProcedure
    \end{algorithmic}
  \end{algorithm}
}

%% file: prol.tex
In the Smoothed Aggregation AMG method, a smoothing step is applied to the tentative prolongation matrix $\mP$ to obtain
the smoothed prolongation matrix $\mP_s$ (see e.g \cite{Vanek_SA_1992})
\[
\mP_s = (\mI-\omega \mD^{-1}\mA)\mP.
\]
Here $\mD$ is the diagonal of $\mA$ and $\omega$ is a dampening parameter.
In order to control the sparsity of $\mP_s$, instead of $\mA$, a filtered, sparser version of $\mA$ could be used here.
As $\mP_s$ still has to preserve kernel vectors, care must be taken in how this is done.
For scalar problems, where the only kernel vector is the constant one, off-diagonal entries can be removed and subtracted from the diagonal entries such that the row-sum remains zero.
For elasticity problems, eliminating off-diagonal entries of $\mA$ without disturbing the rigid body modes is less easy to do.
One such approach can be found in \cite{strong_coupling_SA} (Algorithm 5).
Eliminating off-diagonal entries from the auxiliary matrix $\mhA$, however, is quite straightforward - in order to eliminate $\mhA_{ij}$, we only need to not assemble
the contribution of the edge connecting $v_i$ and $v_j$. This motivates the following definition:
\begin{definition}[Auxiliary Smoothed prolongation]
  \label{def:aux_sp}
  Let $\mP$ be a tentative prolongation as in Definition \ref{def:pwp}.
  Let the set of non-Dirichlet vertices, $\mcV\setminus\mcD$, be split into two sets, $\mathcal{S}$ and $\mathcal{A}$.
  For $v_i\in\mathcal{A}$, define some set $\mathcal{F}^i$ of "filtered neighbors". Then, define the filtered auxiliary matrix by
  \[
    \mhA^0_{ii} :=
    \begin{cases}
      \mA_{ii} &\quad v_i\in\mathcal{S}, \\
      \sum_{l\in\mathcal{F}^i} Q^{li,T}E^{il}Q^{li} &\quad v_i\in\mathcal{A}, \\
      \mI &\quad v_i\in\mcD,
    \end{cases}
    \quad\quad\text{and}\quad\quad
    \mhA^0_{ij} :=
    \begin{cases}
      \mA_{ij} &\quad v_i\in\mathcal{S}, \\
      -Q^{ij,T}M^{ij}Q^{ji} &\quad v_i\in\mathcal{A} \land j\in \mathcal{F}^i,\\
      0 &\quad \text{otherwise}.
    \end{cases}
  \]
  The auxiliary smoothed prolongation is then given as
  \[
    \mP_s := (\mI-\omega (\mhD^0)^{-1}\mhA^0)\mP,
  \]
  where $\mhD^0$ is the diagonal of $\mhA^0$.
\end{definition}

\begin{remark}
  Splitting vertices into $\mathcal{S}$ (standard vertices) and $\mathcal{A}$ (auxiliary vertices) allows for a hybrid approach. 
  For vertices where no filtering is required, e.g., where $\mA$ has only few row entries, the actual matrix can be used, and for other vertices edge matrices are utilized.
\end{remark}

\begin{remark}
  For the elascitity problem, the set $\mathcal{F}^i$ can be chosen such that $\mhA^0_{ii}$ is singular.
  Instead of forcing all $\mathcal{F}^i$ to contain sufficiently many neighbors, it is feasible in practice to use
  $(\mhD^{0})^{\dagger}$ instead of $(\mhD^{0})^{-1}$.
  A final option would be to require every nonempty $\mathcal{F}^i$ to contain at least neighbor $v_j$ that is in the same agglomerate as $v_i$, and to then regularize $M^{ij}$ such that it has full rank which makes $\mhD^0_{ii}$ regular.
\end{remark}

In a distributed parallel setting, the domain, and thus the mesh, is partitioned, with every process having its own assigned subdomain.
Now, usually, as coarse grid matrices become less sparse, they will also have many "off-processor" entries connecting degrees of freedom inside different subdomains.
This way, even processors that do not share a boundary on the finest level can come into contact with one another on coarse levels which increases, in addition to the communication volume, also the amount of messages to be sent.
One approach to combat this effect is to somehow remove entries from coarse grid matrices themselfs, resulting in non-Galerkin methods, \cite{Falgout2014, Treister2015, Bienz2016}.
In contrast, there is a particularly interesting optimization that can be done with the auxiliary smoothed prolongation to address it by construction.

There is a natural half order "$\preceq$" among vertices, where vertex $v_i\preceq v_j$ iff. the set of processors that share $v_j$ is a superset of those that share $v_i$.
We now assume that the coarsening algorithm respects this half order in the sense that for any agglomerate $C\in\mcC$, all containted vertices are comparable.
As we have not touched on the topic of parallel coarsening at all, let it just be said that this means that there are no agglomerates crossing subdomain boundaries
and that this is given in our parallelized version of Algorithm \ref{alg:pwagg}.
Then, imposing the additional restriction $v_j\in\mathcal{F}^i\Rightarrow v_i\preceq v_j$ on filtered neighbor sets entails that transport of information only happens from "more global" vertices to "more local" ones, which motivates the label "hierarchic" for such a prolongation.
A less precise but more visual description would be that we allow coarse grid basis functions associated to boundary degrees of freedom to grow into the interior of attached subdomains, but prevent those in the interior of any of those to touch the boundary.

In this way, a hierarchic prolongation only produces new entries in the coarse grid matrix that connect degrees of freedom on subdomain boundaries with other degrees of freedom contained in either of the attached subdomains.
This guarantees that processes never come into contact with others that are not already neighbors on the fine level.
For the sparse matrix format used in our implementation, where every process stores a diagonal block of the global matrix (with overlapping blocks), the coarse grid matrix can be computed without any communication by simply locally computing the Galerkin projection of the diagonal blocks, using the appropriate sub-block of the prolongation matrix.
Additionally, hybrid smoothers (such as the $\ell_1$ smoother introduced in \cite{UPSM}) can be implemented such that their application only requires the exchange of partial residual entries (which can be computed locally).
In that case, the message size does not depend on the number of nonzero entries in rows belonging to shared degrees of freedom, but only on their number, which is not increased by a hierarchic prolongation.
Therefore, for the combination of a parallel sparse matrix format based on overlapping diagonal blocks and a hybrid smoother, smoothing the prolongation comes at no additional communication cost.

Of course, there is apparently one big disadvantage to this approach.
As the number of levels increases, and thus the number of degrees of freedom per subdomain decreases, the additional restriction we impose becomes very strict to the point that most $\mathcal{F}^i$ will be empty and we are left with the tentative prolongation we started out with.
The solution is to merge and re-distribute subdomains (to fewer processes) periodically, which frees up additional connections to be taken into account for smoothing the tentative prolongation.
This is convenient in practice anyways, as otherwise, on coarse levels, when the number of degrees of freedom per process is very low, the cost of communication can dominate the cost of computation.
There is, however, a balance to be struck here.
When subdomains are merged too aggressively or too frequently, computation costs on coarse levels can dominate the costs on the finest level locally.
Then, even though the method could have a perfectly acceptable theoretical cost, much of the work has to be done on few processors and performance suffers.
Ultimately, however, in our opinion, this difficulty is more than made up for by the performance improvements that can be gained, especially when one considers that this kind of load balancing issue arises on coarse levels anyways.
As a final argument in favor of a hierarchical approach, we want to stress the ease of implementation of the sparse matrix product required to compute coarse grid matrices for which one merely needs an efficient sequential (or shared memory parallel) sparse matrix multiplication for completely local sparse matrices.

Finally, there is one aspect to using auxiliary smoothed prolongations that cannot be addressed with the theory in this paper.
Theorem \ref{thm:cemats} shows how to define coarse grid edge matrices such that $\mA_c=\mP^T\mA\mP \sim \mhA_c = \mP^T\mhA \mP$.
However, using $\mP_s$ instead of the tentative $\mP$, the coarse matrix is $\mA_c = \mP_s^T\mA\mP_s \neq \mP^T\mA\mP$.
This coarse grid matrix is denser than the previous one, and we make no statements about edge matrices that induce a coarse level auxiliary matrix equivalent to the new one.
It would be interesting to further study the problem of finding coarse grid edge matrices in this case.
Possibly, techniques such as used in \cite{CMIII} to find the proper scaling of the finest level edge matrices based on off-diagonal entries could work for that.
While this would potentially improve the quality of smoothed prolongations, the coarsening algorithm might even profit from this distinction since $\mhA_c$ is much sparser than $\mA_c$.
This effect is compounded for many levels.
The sparsity of the auxiliary matrices remains approximately the same, while the density of coarse grid matrices keeps growing due to repeated Galerkin projection with (denser than tentative) auxiliary smoothed prolongations.
This means that the coarsening algorithm only considers ever matches coarse grid vertices that are actually "geometrically" next to one another.

%% file: res.tex

We now present computational results in an attempt to demonstrate both the robustness of the coarsening algorithm presented in Section \ref{sec:crs} as well as the scalability of out implementation.
In all cases, the Conjugate Gradient method with one V-cycle of the auxiliary topology AMG method as a preconditioner was used to reduce the error by a factor of $10^6$.
The smoother used is a block version of the $\ell_1$ smoother introduced in \cite{UPSM}, which was chosen for its parallelizability.
Blocks are made up of all degrees of freedom associated with a single vertex.
We consider three problems, with the first two being of a simple nature, intended to show scalability of the implementation for both the scalar and vector case.
The last problem is the three dimensional analogon of the problem used in Figure \ref{fig:2dc_rob} and is as such explicitly constructed to make traditional coarsening algorithms fail.
It is intended to demonstrate robustness of the coarsening algorithm even for difficult problems on a large scale.

The AMG method has been implemented as an extension library to the C++ open source mesh generator and finite element software Netgen/NGSolve \cite{Netgen, cpp_ngs}, see also the webpage of that project, \cite{ngs_webpage}.
It is available as research code through the "NgsAMG"-project on GitHub, \cite{ngs_amg_github}.
The computations were performed on the Vienna Scientific Cluster in its newest 2020 iteration (VSC4).
For all cases, a series of unstructured tetrahedral meshes of varying size have been generated for the same geometry, and the number of vertices per processor has been kept as constant as feasible.
Due to technical limitations, the larger meshes could not be generated beforehand, instead a coarser mesh had to be generated, distributed and then refined in parallel at execution time.
For these reasons we were unable to keep the number of vertices per processor perfectly constant, however it generally lies within 5\% variation off the goal value.

\paragraph{Scalar model problem} 
First, in order to demonstrate scalability of the implementation for the simplest scalar case, we consider Poisson's equation, Problem \ref{prob:gg}, on the unit cube $\Omega=[0,1]^3$ with $\alpha=1$ and Dirichlet conditions on the entire boundary.
The coarsening algorithm uses four passes on the first two levels and three passes after that and uses the scalar criterion $\muc$ only.
Coarsening is stopped when the number of vertices reaches the mininmum of a static threshold ($1600$) and the initial number of vertices divided by a factor ($1250$).
The number of entries per row of the first three auxiliary smoothed prolongation is capped by four.
For later ones, the cap is four for vertices where edge matrices are used, and six (for the fourth and fifth) or eight (any later ones) for those where real matrix entries are used.
The results are listed in Table \ref{tab:ucu}.
For each computation we give the number of nodes (\#N), processors (\#P) and degrees of freedom (as a multiple of the number of vertices, \#V), as well as the number of vertices per processor (\#V/\#P).
The number of AMG levels (\#L) and the number of vertices on the coarsest one ($n_c$) and both the grid, or vertex complexity (VC), the sum of the number of vertices on all levels divided by the initial number of vertices and the operator complexity (OC), the number of nonzero entries of matrices on all levels divided by those of the finest level matrix describe the AMG iteration itself.
Finally, the number of CG iterations necessary to reduce the error by a factor of $10^6$ (\#its), as well as the time for the setup of the AMG levels (tsup) and the
CG iterations (tsol) are complemented by $\text{eff}_{\text{sup}}$ and $\text{eff}_{\text{sol}}$, which is the sequential time divided by either tsup or tsol.

\begin{table}[!ht]
  \centering
  \input{tab_ucu.tex}
  \caption{Results for the scalar model problem. }
  \label{tab:ucu}
\end{table}


\paragraph{Vector model problem} 
Moving on to vector problems, we consider Problem \ref{prob:epseps} on a beam $\Omega=[0, 10]\times[0,1]^2$ that is fixed on the left side  $\Gamma_D = \{0\}\times[0,1]\times[0,1]$. Suitable for a model problem, we choose constant parameters $\mu=1$ and $\lambda=0$.
This time we properly make use of all three measures $\mua, \mub$ and $\muc$, use five passes on the first and four on subsequent levels, and set the minimum threshold for the coarsest level to the minimum of $1000$ and the initial number of vertices divided by $10^4$.
The entries per row of the prolongations is capped by four on all levels for vertices where edge matrices are used and eight (the first level) or six (subsequent levels) for others.
Note that, with a comparable number of degrees of freedom per process to that of the sclar model problem, although the matrix is denser, and although the more expensive
criteria $\mub$ and $\muc$ are used for coarsening, the setup time is in the same order of magnitude.

{
  \begin{table}[!ht]
    \centering
    \input{tab_3dbeam.tex}
    \caption{Results for the vector model problem on a three dimensional beam.}
    \label{tab:3dbeam}
  \end{table}
}

\begin{center}
  \begin{figure}[!ht]
    \centering
    \includegraphics[scale=.8]{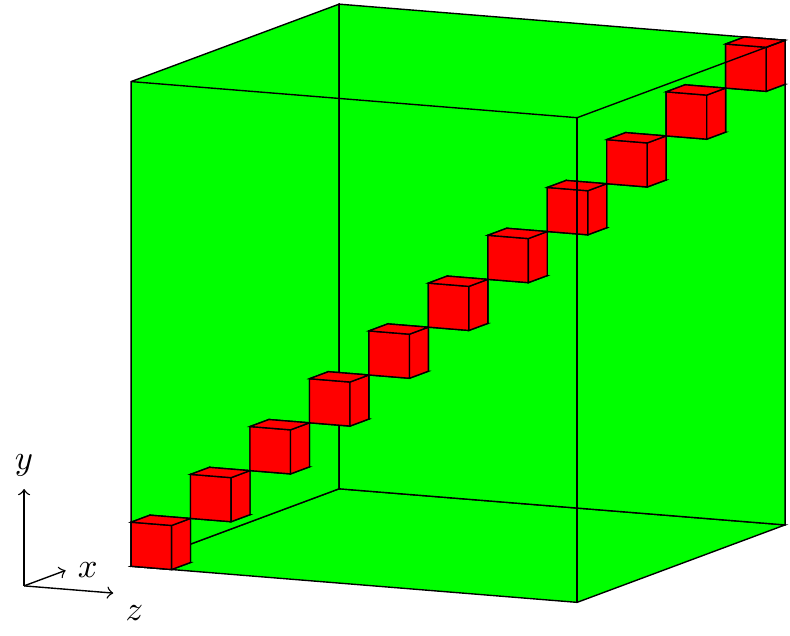}
    \caption{The geometry for the difficult vector problem, $\Omega_1$ in green and $\Omega_2$ in red.}
    \label{fig:crs_fail}
\end{figure}
  \end{center}

\paragraph{A difficult vector problem}
The final problem is the three dimensional analogon of the problem used in Figures \ref{fig:2dc_rob} and \ref{fig:2dc_norob} to illustrate the necessity
for the robust criteria $\mua$ and $\mub$.
It is as such specifically constructed to defeat traditional, scalar, criteria and features a large jump in coefficients $\mu, \lambda$.
The domain $\Omega=[0,1]^3$ is split into a subdomain $\Omega_1$ and its complement $\Omega_2$. As depicted in Figure \ref{fig:crs_fail},
$\Omega_2$ is the union of eleven "boxes" that touch in their corners, 
$$
\Omega_2 = \bigcup_{i=1}^{11} \left[\frac{i-1}{11}, \frac{i}{11}\right]^3.
$$
With $\mu=\lambda$ we are still far away from the almost incompressible case, however the large jump between $\mu=1$ in $\Omega_1$ and $\mu=10^4$ in $\Omega_2$ is problematic.
The difficulty comes from the fact that the coarsening algorithm has to correctly identify that no agglomerates spanning multiple of these boxes must be formed.
The exception is that, once all vertices within two boxes on the right end of the chain are contained within one agglomerate each, they can be combined.
On the next level, another single box can be added, etc, until all but the one first one, that touches the Dirichlet boundary, are agglomerated.
That one is different in that its connection to the Dirichlet boundary $\Gamma_D=\{0\}\times[0,1]^2$, once all (non-Dirichlet) vertices in it are agglomerated together, is represented by an edge weight matrix that has six large eigenvalues, while edge weight matrices for any connection between two such box agglomerates have three large and three small eigenvalues.
Only the robust criteria $\mua$ and $\mub$ can correctly identify this.
Parameters for the AMG method are the same as for the vector model problem, except that, as this is a hard problem, only four passes were done on the first level and three afterwards.
Results are listed in Table \ref{tab:co3d}.

{
  \begin{table}[!ht]
    \centering
    \input{tab_co3d.tex}
    \caption{Results for the difficult vector problem.}
    \label{tab:co3d}
  \end{table}
}

%% file: tab_ucu.tex
\begin{tabular}{l|l|l|l|l|l|l|l|l|l|l|l|l}
\hline\hline
 $\#N$ & $\#P$ & $\#V$ & $\#V/\#P$ & $\#L$ & $n_c$ & VC & OC & $\#\text{its}$ & tsup & $\text{eff}_{\text{sup}}$ & tsol & $\text{eff}_{\text{sol}}$ \\
\hline\hline
1 & 1 & 151.64K & 151.64K & 4 & 60 & 1.06 & 1.20 & 13 & 4.35 & 1.00 & 0.24 & 1.00 \\
1 & 4 & 600.52K & 150.13K & 4 & 267 & 1.06 & 1.22 & 14 & 5.16 & 0.84 & 0.31 & 0.77 \\
1 & 19 & 2.81M & 148.05K & 4 & 1362 & 1.06 & 1.24 & 15 & 6.21 & 0.70 & 0.40 & 0.60 \\
2 & 49 & 7.08M & 144.42K & 5 & 448 & 1.07 & 1.25 & 19 & 6.81 & 0.64 & 0.63 & 0.38 \\
3 & 138 & 20.78M & 150.59K & 5 & 1427 & 1.07 & 1.25 & 18 & 7.63 & 0.57 & 1.26 & 0.19 \\
8 & 362 & 53.87M & 148.82K & 6 & 478 & 1.07 & 1.25 & 26 & 7.39 & 0.59 & 1.67 & 0.14 \\
12 & 537 & 69.57M & 129.55K & 6 & 645 & 1.07 & 1.25 & 26 & 8.80 & 0.49 & 1.65 & 0.15 \\
24 & 1.108 & 164.85M & 148.78K & 6 & 1543 & 1.07 & 1.26 & 26 & 8.75 & 0.50 & 1.95 & 0.12 \\
43 & 2.049 & 296.44M & 144.68K & 7 & 377 & 1.07 & 1.25 & 32 & 8.05 & 0.54 & 2.11 & 0.11 \\
65 & 3.072 & 446.47M & 145.33K & 7 & 592 & 1.07 & 1.26 & 33 & 8.30 & 0.52 & 2.28 & 0.11 \\
84 & 4.007 & 581.6M & 145.15K & 7 & 771 & 1.07 & 1.26 & 34 & 8.53 & 0.51 & 2.35 & 0.10 \\
185 & 8.866 & 1.31B & 148.2K & 8 & 235 & 1.07 & 1.26 & 37 & 8.80 & 0.49 & 2.51 & 0.10 \\
219 & 10.505 & 1.41B & 133.79K & 8 & 372 & 1.07 & 1.25 & 37 & 17.12 & 0.25 & 2.41 & 0.10 
\end{tabular}

%% file: tab_3dbeam.tex
\begin{tabular}{l|l|l|l|l|l|l|l|l|l|l|l|l}
\hline\hline
 $\#N$ & $\#P$ & $3\cdot\#V$ & $\#V/\#P$ & $\#L$ & $n_c$ & VC & OC & $\#\text{its}$ & tsup & $\text{eff}_{\text{sup}}$ & tsol & $\text{eff}_{\text{sol}}$ \\
\hline\hline
1 & 1 & 128.53K & 42.84K & 4 & 22 & 1.08 & 1.26 & 20 & 15.95 & 1.00 & 0.68 & 1.00 \\
1 & 3 & 431.37K & 47.93K & 4 & 32 & 1.05 & 1.17 & 16 & 16.31 & 0.98 & 0.63 & 1.08 \\
1 & 11 & 1.35M & 41.02K & 4 & 111 & 1.05 & 1.18 & 16 & 14.60 & 1.09 & 0.61 & 1.11 \\
1 & 19 & 2.3M & 40.27K & 4 & 198 & 1.05 & 1.20 & 17 & 15.46 & 1.03 & 0.75 & 0.91 \\
1 & 28 & 3.27M & 38.97K & 4 & 283 & 1.06 & 1.23 & 17 & 16.31 & 0.98 & 0.98 & 0.69 \\
1 & 37 & 4.24M & 38.15K & 4 & 383 & 1.06 & 1.23 & 17 & 16.67 & 0.96 & 1.30 & 0.52 \\
2 & 90 & 10.4M & 38.5K & 4 & 921 & 1.06 & 1.23 & 17 & 18.71 & 0.85 & 1.65 & 0.41 \\
4 & 147 & 17.73M & 40.2K & 5 & 131 & 1.06 & 1.22 & 21 & 19.08 & 0.84 & 1.79 & 0.38 \\
5 & 218 & 25.46M & 38.92K & 5 & 195 & 1.06 & 1.27 & 24 & 21.02 & 0.76 & 2.53 & 0.27 \\
25 & 1.159 & 136.51M & 39.26K & 6 & 151 & 1.06 & 1.27 & 33 & 24.05 & 0.66 & 3.93 & 0.17 \\
28 & 1.309 & 145.87M & 37.14K & 5 & 995 & 1.06 & 1.22 & 22 & 24.23 & 0.66 & 2.57 & 0.26 \\
36 & 1.697 & 200.69M & 39.42K & 7 & 135 & 1.07 & 1.28 & 46 & 25.93 & 0.62 & 5.78 & 0.12 \\
48 & 2.258 & 259.88M & 38.36K & 6 & 160 & 1.06 & 1.25 & 30 & 25.76 & 0.62 & 3.82 & 0.18 \\
116 & 5.544 & 647.46M & 38.93K & 6 & 416 & 1.06 & 1.23 & 34 & 25.43 & 0.63 & 4.15 & 0.16 \\
210 & 10.039 & 1.14B & 37.83K & 7 & 116 & 1.06 & 1.21 & 40 & 29.33 & 0.54 & 5.07 & 0.13 \\
307 & 14.724 & 1.59B & 36.08K & 7 & 249 & 1.05 & 1.19 & 39 & 30.47 & 0.52 & 4.38 & 0.16 
\end{tabular}

%% file: tab_co3d.tex
\begin{tabular}{l|l|l|l|l|l|l|l|l|l|l|l|l}
\hline\hline
 $\#N$ & $\#P$ & $3\cdot\#V$ & $\#V/\#P$ & $\#L$ & $n_c$ & VC & OC & $\#\text{its}$ & tsup & $\text{eff}_{\text{sup}}$ & tsol & $\text{eff}_{\text{sol}}$ \\
\hline\hline
1 & 1 & 87.73K & 29.25K & 5 & 22 & 1.11 & 1.37 & 17 & 9.06 & 1.00 & 0.66 & 1.00 \\
1 & 4 & 448.9K & 37.41K & 5 & 46 & 1.10 & 1.39 & 16 & 13.41 & 0.68 & 1.04 & 0.63 \\
1 & 9 & 1.04M & 38.48K & 5 & 93 & 1.10 & 1.39 & 16 & 15.06 & 0.60 & 1.23 & 0.54 \\
2 & 51 & 5.24M & 34.22K & 5 & 382 & 1.10 & 1.39 & 21 & 17.46 & 0.52 & 2.40 & 0.28 \\
2 & 69 & 8.1M & 39.13K & 6 & 100 & 1.11 & 1.44 & 23 & 18.20 & 0.50 & 3.38 & 0.20 \\
3 & 115 & 13.55M & 39.29K & 6 & 191 & 1.11 & 1.47 & 24 & 20.49 & 0.44 & 3.50 & 0.19 \\
4 & 162 & 14.85M & 30.57K & 6 & 192 & 1.09 & 1.36 & 29 & 24.15 & 0.38 & 5.29 & 0.12 \\
6 & 250 & 28.32M & 37.76K & 6 & 307 & 1.11 & 1.42 & 29 & 27.76 & 0.33 & 5.76 & 0.11 \\
12 & 529 & 62.19M & 39.19K & 7 & 195 & 1.12 & 1.50 & 29 & 20.85 & 0.43 & 5.37 & 0.12 \\
22 & 1.029 & 116.31M & 37.68K & 7 & 173 & 1.10 & 1.43 & 27 & 23.09 & 0.39 & 6.32 & 0.10 \\
42 & 2.005 & 225.6M & 37.51K & 7 & 370 & 1.10 & 1.42 & 41 & 48.51 & 0.19 & 11.96 & 0.06 \\
63 & 2.994 & 329.21M & 36.65K & 8 & 108 & 1.10 & 1.40 & 45 & 49.09 & 0.18 & 13.01 & 0.05 \\
155 & 7.419 & 850.97M & 38.23K & 8 & 243 & 1.10 & 1.41 & 37 & 34.32 & 0.26 & 9.71 & 0.07 \\
226 & 10.831 & 1.27B & 38.93K & 8 & 359 & 1.10 & 1.41 & 34 & 31.17 & 0.29 & 9.69 & 0.07 
\end{tabular}

%% file: preprint_elamg.bbl
\begin{thebibliography}{10}

\bibitem{ngs_webpage}
Netgen/ngsolve software.

\bibitem{ngs_amg_github}
Ngsamg software.

\bibitem{UPSM}
Allison Baker, Robert D.~Falgout, Tzanio Kolev, and Ulrike Yang.
\newblock Multigrid smoothers for ultraparallel computing.
\newblock {\em {SIAM} J Sci Comput}, 33:2864--2887, 2011.

\bibitem{Bienz2016}
Amanda Bienz, Robert~D. Falgout, William Gropp, Luke~N. Olson, and Jacob~B.
  Schroder.
\newblock Reducing parallel communication in algebraic multigrid through
  sparsification.
\newblock {\em {SIAM} J Sci Comput}, 38(5):S332--S357, January 2016.

\bibitem{Braess2007}
Dietrich Braess.
\newblock {\em Finite Elements: Theory, Fast Solvers, and Applications in Solid
  Mechanics}.
\newblock Cambridge University Press, 3 edition, 2007.

\bibitem{Brandt1982AlgebraicM}
A.~Brandt, S.~McCormick, and J.~Ruge.
\newblock Algebraic multigrid (amg) for automatic multigrid solutions with
  application to geodetic computation.
\newblock Technical report, Inst. for Computational Studies, Fort Collins, CO,
  1982.

\bibitem{Brezina_eAMG_2001}
M.~Brezina, A.~J. Cleary, R.~D. Falgout, V.~E. Henson, J.~E. Jones, T.~A.
  Manteuffel, S.~F. McCormick, and J.~W. Ruge.
\newblock Algebraic multigrid based on element interpolation ({AMGe}).
\newblock {\em {SIAM} J Sci Comput}, 22(5):1570--1592, 2001.

\bibitem{strong_coupling_SA}
Tony~F. Chan and Petr Van{\v{e}}k.
\newblock Detection of strong coupling in algebraic multigrid solvers.
\newblock In {\em Lecture Notes in Computational Science and Engineering},
  pages 11--23. Springer Berlin Heidelberg, 2000.

\bibitem{Falgout2004}
R.~D. Falgout and P.~S. Vassilevski.
\newblock On generalizing the algebraic multigrid framework.
\newblock {\em {SIAM} J Numer Anal}, 42(4):1669--1693, January 2004.

\bibitem{Falgout2014}
Robert~D. Falgout and Jacob~B. Schroder.
\newblock Non-galerkin coarse grids for algebraic multigrid.
\newblock {\em {SIAM} J Sci Comput}, 36(3):C309--C334, January 2014.

\bibitem{Sch_elpre_2001}
G~Haase, U.~Langer, S.~Reitzinger, and J.~Sch\"{o}berl.
\newblock Algebraic multigrid methods based on element preconditioning.
\newblock {\em Int J of Comput Math}, 78(4):575--598, 2001.

\bibitem{ef_eamg}
Van~Emden Henson and Panayot~S. Vassilevski.
\newblock Element-free amge: General algorithms for computing interpolation
  weights in amg.
\newblock {\em {SIAM} J Sci Comput}, 23(2):629--650, 2001.

\bibitem{CMIII}
E.~Karer and J.~K. Kraus.
\newblock Algebraic multigrid for finite element elasticity equations:
  Determination of nodal dependence via edge-matrices and two-level
  convergence.
\newblock {\em Int J Numer Meth Eng}, 83(5):642--670, 2010.

\bibitem{CMII}
J.~Kraus.
\newblock Algebraic multigrid based on computational molecules, 2: Linear
  elasticity problems.
\newblock {\em {SIAM} J Sci Comput}, 30(1):505--524, 2008.

\bibitem{CMI}
J.~K. Kraus and J.~Schicho.
\newblock Algebraic multigrid based on computational molecules, 1: Scalar
  elliptic problems.
\newblock {\em Computing}, 77(1):57--75, 2006.

\bibitem{Napov2010}
Artem Napov and Yvan Notay.
\newblock Algebraic analysis of aggregation-based multigrid.
\newblock {\em Numer Linear Algebr}, 18(3):539--564, 2010.

\bibitem{Napov2012}
Artem Napov and Yvan Notay.
\newblock An algebraic multigrid method with guaranteed convergence rate.
\newblock {\em {SIAM} Journal on Scientific Computing}, 34(2):A1079--A1109,
  2012.

\bibitem{RS_1984}
J.~Ruge and K.~St{\"u}ben.
\newblock {\em Efficient Solution of finite difference and finite element
  equations by algebraic multigrid AMG}.
\newblock Arbeitspapiere der GMD. Gesellschaft f. Mathematik u.
  Datenverarbeitung, 1984.

\bibitem{Ruge1987}
J.~W. Ruge and K.~St\"{u}ben.
\newblock 4. algebraic multigrid.
\newblock In {\em Multigrid Methods}, pages 73--130. Society for Industrial and
  Applied Mathematics, January 1987.

\bibitem{Netgen}
Joachim Sch\"{o}berl.
\newblock {NETGEN} an advancing front 2d/3d-mesh generator based on abstract
  rules.
\newblock {\em Comput Vis Sci}, 1(1):41--52, 1997.

\bibitem{cpp_ngs}
Joachim Sch\"{o}berl.
\newblock C++11 implementation of finite elements in ngsolve, 2014.

\bibitem{Treister2015}
Eran Treister and Irad Yavneh.
\newblock Non-galerkin multigrid based on sparsified smoothed aggregation.
\newblock {\em {SIAM} J Sci Comput}, 37(1):A30--A54, January 2015.

\bibitem{Vanek_SA_1996}
P.~Van{\v{e}}k, J.~Mandel, and M.~Brezina.
\newblock Algebraic multigrid by smoothed aggregation for second and fourth
  order elliptic problems.
\newblock {\em Computing}, 56(3):179--196, 1996.

\bibitem{Vanek_SA_1992}
Petr Van{\v{e}}k.
\newblock Acceleration of convergence of a two-level algorithm by smoothing
  transfer operators.
\newblock {\em Appl Math-Czech}, 37(4):265--274, 1992.

\bibitem{Vanek_SA_2001}
Petr Van{\v{e}}k, Marian Brezina, and Jan Mandel.
\newblock Convergence of algebraic multigrid based on smoothed aggregation.
\newblock {\em Numer Math}, 88(3):559--579, 2001.

\bibitem{Vanek_SA_1999}
Petr Van{\v{e}}k, Marian Brezina, and Radek Tezaur.
\newblock Two-grid method for linear elasticity on unstructured meshes.
\newblock {\em {SIAM} J Sci Comput}, 21(3):900--923, 1999.

\bibitem{mlfvass}
Panayot~S. Vassilevski.
\newblock {\em Multilevel Block Factorization Preconditioners}.
\newblock Springer New York, 2008.

\bibitem{amg_ov_Xu2017}
Jinchao Xu and Ludmil Zikatanov.
\newblock Algebraic multigrid methods.
\newblock {\em Acta Numer}, 26:591--721, 2017.

\end{thebibliography}
